\documentclass[10pt,reqno]{amsart}
\usepackage{mathrsfs}
\usepackage{graphicx,epsfig,amsmath,amsfonts,dsfont}
\usepackage{color}
\usepackage{latexsym,amssymb}
\usepackage{lineno,version}
\usepackage{multirow}
\usepackage{bm}
\usepackage{float}
\usepackage{booktabs}
\usepackage{hyperref} 
\usepackage{geometry}
\allowdisplaybreaks

\catcode`\@=11
\theoremstyle{plain}
\@addtoreset{equation}{section}   

\@addtoreset{figure}{section}
\renewcommand\thefigure{\thesection.\@arabic\c@figure}

\newtheorem{thm}{\noindent Theorem}[section]
\newtheorem{rem}{\noindent Remark}[section]
\newtheorem{alg}{\noindent Algorithm}[section]
\newtheorem{prop}{\noindent Proposition}[section]
\newtheorem{lem}{\noindent Lemma}[section]

\newcommand{\ba}{\begin{array}}\newcommand{\ea}{\end{array}}
\newcommand{\be}{\begin{eqnarray}}\newcommand{\ee}{\end{eqnarray}}
\newcommand{\beq}{\begin{equation}}\newcommand{\eeq}{\end{equation}}
\newcommand{\bex}{\begin{eqnarray*}}
\newcommand{\eex}{\end{eqnarray*}}

\def\dps{\displaystyle}

\font\tenbi=cmmib10   at 11 pt
\font\sevenbi=cmmib10 at 9pt
\font\fivebi=cmmib7 at 6pt
\newfam\bifam
\textfont\bifam=\tenbi \scriptfont\bifam=\sevenbi
\scriptscriptfont\bifam=\fivebi

\font\tendb=msbm10 at 12 pt
\font\sevendb=msbm7
\newfam\dbfam
\textfont\dbfam=\tendb \scriptfont\dbfam=\sevendb

\def\div{\mbox{div}}

\begin{document}

\title[Inverse source problem for TFDE with $L^2$-TV regularization]
{Identifying source term in the subdiffusion equation with $L^2$-TV regularization$^*$}
\author[B. Fan \& C.J. Xu]
{Bin Fan$^{1}$
\quad
Chuanju Xu$^{1,\dag}$}
\thanks{\hskip -12pt
${}^*$This research is partially supported by NSFC grant 11971408 and NNW2018-ZT4A06 project.\\
${}^{1}$School of Mathematical Sciences and
Fujian Provincial Key Laboratory of Mathematical Modeling and High Performance
Scientific Computing, Xiamen
University, 361005 Xiamen, China.\\
Email: binfan@stu.xmu.edu.cn (B. Fan), cjxu@xmu.edu.cn (C.J. Xu)\\
${}^{\dag}$Corresponding author}



\maketitle

\begin{abstract}
In this paper, we consider the inverse source problem for the time-fractional diffusion equation,
which has been known to be an ill-posed problem.
To deal with the ill-posedness of the problem, we propose to transform the problem
into a regularized problem with $L^2$ and total variational (TV) regularization terms.
Differing from the classical Tikhonov regularization with $L^2$ penalty terms,
the TV regularization is beneficial for reconstructing discontinuous or
piecewise constant solutions. The regularized problem is then approximated
by a fully discrete scheme.
Our theoretical results include: estimate of
the error order between the discrete problem and the continuous direct problem;
the convergence rate of the discrete regularized solution to the target source term; and
the convergence of the regularized solution with respect to the noise level.
Then we propose an accelerated primal-dual iterative algorithm
based on an equivalent saddle-point reformulation of the discrete regularized model.
Finally, a series of numerical tests are carried out to
demonstrate the efficiency and accuracy of the algorithm.

\end{abstract}

\section{Introduction}\label{sec:1}

The fractional diffusion equations, or more general fractional partial differential equations,
have been the subjects of many research, and received
increasing attention due to its wide range of applications in science and
engineering. Unlike ordinary derivatives, fractional derivatives are non-local in nature and
are capable of modeling memory or nonlocal effects.
Various models using fractional derivatives have been proposed
and there has been significant interest in developing numerical methods for their solutions;
see, e.g., monographs \cite{Podlubny-1999,Mainardi-2010,diethelm2010analysis,Handbook2019}
and the references therein.

In this paper,
we are interested in the time-fractional diffusion equation (TFDE). More precisely,
our point of interest is to solve the inverse problem of identifying the source $G(x,t)$
in the following equation:
\begin{eqnarray}
   && \left\{\begin{array}{ll}
               \partial_t^\alpha u(x,t)+L u(x,t)=G(x,t), & x\in\Omega,\quad t\in(0,T], \\
               u(x,0)=u_0(x),& x\in\Omega, \\
               u(x,t)=0, &  x\in\partial\Omega,\quad t\in[0,T],
             \end{array}\right. \label{general-dir-prob}
\end{eqnarray}
where $\Omega$ is a bounded and convex polygonal domain in $\mathds{R}^d$ ($d=1,2$),
$T>0$ is a fixed time,
$\partial_t^\alpha$ denotes the Caputo fractional derivative of order $\alpha$ defined by
\begin{eqnarray}
 \partial_t^\alpha u(x,t):=\frac{1}{\Gamma(1-\alpha)}\int_0^t \frac{\partial u(x,s)}{\partial s} \frac{\mathrm{d}s}{(t-s)^{\alpha}},\quad 0<\alpha<1,\quad 0<t\leq T
\label{det-caputo-1}
\end{eqnarray}
with $\Gamma(\cdot)$ being the Gamma function.
We suppose $L$ is symmetric uniformly elliptic operator,
which can be $-\Delta$ or the following more general operator:
\begin{eqnarray*}
 Lv(x,t)= - \sum\limits_{i=1}^d\frac{\partial}{\partial x_i}\left(\sum\limits_{i=1}^da_{ij}(x)\frac{\partial}{\partial x_j}v(x,t) \right)+b(x)u(x,t)
\end{eqnarray*}
with the coefficients $\{a_{ij}\}$ and $b$ satisfying
\begin{eqnarray*}
   &&  a_{ij}=a_{ji},\quad a_{ij}\in C^1(\bar{\Omega}),\quad 1\leq i,j\leq d,\\
   && a_0\sum\limits_{i=1}^d\xi_i^2\leq \sum\limits_{i,j=1}^da_{ij}(x)\xi_i\xi_j,\quad \forall x\in\bar{\Omega},
   \quad \forall \xi\in\mathds{R}^d,\quad a_0>0,\\
   && b\in C(\bar{\Omega}),\quad b(x)\ge 0,\quad \forall x\in\bar{\Omega}.
\end{eqnarray*}

If the functions $u_0(x)$ and $G(x,t)$ are given appropriately, the problem (\ref{general-dir-prob}) is a direct problem.
The direct problems for time fractional diffusion equations have been studied extensively; see, e.g., recent papers
\cite{ZX19,Jin-2019} and the references therein.
However, the research on the related inverse problems is relatively sparse.
As it is well know, these inverse problems are generally ill-posed, similar to the inverse problems associated to
integer order diffusion equations \cite{Engl-Hanke-Beyvayer-1996}.
By using eigenfunction expansion of elliptic operator in space, Liu and Yamamoto \cite{Liu-Yam-2010} proposed a quasi-reversibility method for an initial value problem of \eqref{general-dir-prob}.
Wang et al. \cite{Wang-Wei-2013} transformed the initial value inverse problem
into a Fredholm integral equation of the first kind,
then solved the transformed problem by Tikhonov regularization method.
In \cite{Wei-Wang-2014}, Wei and Wang considered a space-dependent source term in the time-fractional diffusion equation with a noisy final time data.
They proposed a modified quasi-boundary value regularization method to deal with the inverse source problem
and obtained the convergence rates
under a priori and a posteriori regularization parameter choice rule.
A common feature of these methods is that they require the spectral information of the operator $-L$.
This limits the applicability of the method to a larger extent, especially to general domain $\Omega$, in which
exact spectral information is usually unknown.
In recent years, a number of variational methods have also been developed. The variational methods
consist in approaching the solutions of problems based on weak formulations.
Ye and Xu \cite{Ye-Xu-2013} transformed an initial value problem into an optimal control problem based on Tikhonov functional and proposed a space-time spectral method. The error estimate was given
assuming sufficiently smooth state and control variables.
Similarly, under the optimal control framework, the finite element approximation of a class of inverse source problems
was proposed by Zhou and Gong \cite{Zhou-Gong-2016} and Jin et al \cite{Jin-Li-2020}.
In the absence of control constraints, Jiang et al. \cite{Jiang-Liu-2020} carried out a convergence analysis
for the numerical solution
of the spatial component in the source term of the time-fractional diffusion equation.
Wei et al. \cite{Wei-Xian-2019,Yan-Wei-2019} also formulated the initial value problem into a variational problem by Tikhonov functional, and obtained an approximation to the minimizer by using a conjugate gradient method.

It is known that the classical Tikhonov methods with standard penalty terms such as $\|\cdot\|_{L^2(\Omega)}^2$ suffer from the over-smoothing phenomena.
In image processing, a common way to overcome the over-smoothing is to use
total variation (TV) regularization,
which is advantageous to recover non-smooth or discontinuous solutions.
However, compared with the wide application of TV regularization method in image processing \cite{Chambolle-2004,Chan-Mar-2000},
its use in inverse problems of differential equations is still sparse;
see, e.g., \cite{Wang-Liu-2013,Hinze-Quyen-2019,Hinze-Kalt-2018}.
One of the reasons for the limit application of TV regularization method in PDEs-based inverse problems
is mainly due to the non-differentiability of the TV-term,
which brings difficulties to theoretical analysis and numerical calculation.
In the present work, we will make attempt to construct a TV regularization method for the source inverse
problem mentioned above. Our main contribution is as follows:
1) First, the underlying problem is discretized by using the standard piecewise linear finite element in space
and a finite difference scheme in time, then we transform the discrete inverse problem
into an optimization problem with $L^2$-TV regularization.
2) Under appropriate assumptions, we establish some convergence results for the proposed method.
Unlike many of the previous work, our analysis is carried out in the finite dimensional space.
3) An accelerated linearized primal-dual algorithm is proposed to solve the discrete minimization problem.
The benefit of our work is that the proposed method does not require explicit spectrum information
of the operator $L$, thus is extendable to problems in general domains.

The rest of this paper is organized as follows. In Section \ref{sec:2}, we review the some existing results of direct problem and the useful properties of the bounded total variation spaces, and transform the original problem into an optimization problem with $L^2$-TV regularization.
Some useful approximation results for $BV$ spaces are given in \ref{sec:x}.
In Section \ref{sec:3}, we discretize the continuous optimization problem by the finite element method, and give the main convergence results of this paper. In Section \ref{sec:4}, we transform the discrete optimization problem into an equivalent saddle point problem, and analyze the optimality condition. On this basis, we give an effective algorithm for solving saddle point problem.
Numerical implementation and some results are reported in
Section \ref{sec:5}. Finally, the conclusions are given in Section \ref{sec:6}

\section{Problem and preliminaries}\label{sec:2}

Without loss of generality, we consider the equation \eqref{general-dir-prob} with $L=-\Delta$. As an inverse source
problem, our focus is recovering the source term in form of separated variables. Precisely,
we consider the time-fractional diffusion equation
\be
\partial_t^\alpha u(x,t)-\Delta u(x,t)=\mu(t)f(x), & x\in\Omega,\quad t\in(0,T],
\label{dir-prob}
\ee
subject to the homogeneous initial boundary conditions,
where the spatial component $f(x)$ of the source term models the spatial distribution,
and the time component $\mu(t)$ describes the time evolution pattern.
The inverse problem that we are concerned with reads:
 (\textbf{IP}) Given $\mu(t)$ and the final state observation $g(x)$,
 determine $f(x)$ such that
 \begin{eqnarray}
 \|g(x)-u(x,T)\|_{L^2(\Omega)}\leq\delta, \label{level-noisy-g}
\end{eqnarray}
where $u(x,T)$ is the solution of \eqref{dir-prob} associated to the source $\mu(t)f(x)$,
$\delta>0$ is a noise level which may be presented in the observation data $g$.

Throughout the paper we use $c$, with or without subscripts or bars, to mean generic
positive constants, which may not be the same at different occurrences.
We use the standard notations of Sobolev spaces $W^{k,p}(\Omega)$,
$W_0^{k,p}(\Omega), k\ge 0, p\ge 1$, and the corresponding norms or seminorms $\|\cdot\|_{W^{k,p}(\Omega)}$,
$|\cdot|_{W^{k,p}(\Omega)}$, etc \cite{Ambrosio-F-P-2000,Clarlet-1978}. In particular,
we denote $W^{0,p}(\Omega)$ by $L^p(\Omega)$, $\|\cdot\|_{W^{0,p}(\Omega)}$ by $\|\cdot\|_{L^p}$.
The norm of a $\mathds{R}^d$-valued function $\rho\in L^p(\Omega;\mathds{R}^d)$ is defined by $\|\rho\|_{L^p(\Omega)}:=\left(\sum\limits_{i=1}^d\|\rho_i\|_{L^p(\Omega)}^p \right)^{1/p}$,
where $\rho_i\in L^p(\Omega)$ is the $i$-th component of $\rho$.
The vector containing all weak partial derivatives of order $k$ is denoted by $D^kf$. In particular,
$D^1f$, the weak gradient of $f$, is also denoted by $\nabla f$.
With these notations, we have $|f|_{W^{k,p}}=\|D^k f\|_{L^p(\Omega)}$ for any $f\in W^{k,p}(\Omega)$.

\subsection{Reformulation of original problem and ill-posedness}

First we recall the representation of the solution to Problem (\ref{dir-prob}), which will be used in the analysis.
Let $\{\lambda_j\}_{j=1}^{\infty}$ and $\{\varphi_j\in H^2(\Omega)\cap H_0^1(\Omega)\}_{j=1}^{\infty}$
are respectively the eigenvalues and associated eigenfunctions of
the operator $-\Delta$ subject to the homogeneous Dirichlet boundary condition, i.e.,
\begin{eqnarray}
  -\Delta \varphi_j=\lambda_j\varphi_j\ \ \mbox{in } \Omega,\quad\text{and}\quad \varphi_j\big|_{\partial\Omega}=0.
  \label{varphi-n}
\end{eqnarray}
It is well known that the eigenvalues $\{\lambda_j\}$ satisfy
\begin{eqnarray}
 0<\lambda_1\leq\lambda_2\leq\cdots\leq\lambda_j \to +\infty \mbox{ as } j\to\infty,
   \label{lambda-n}
\end{eqnarray}
and the eigenfunctions $\{\varphi_j\}$ form an orthonormal basis of $L^2(\Omega)$.
By using the Laplace transform and convolution rule, the solution of (\ref{dir-prob}) can be formally represented by \cite{Sak-Yam-2011,Jin-2019}
\begin{eqnarray}\label{dir-solution}
  u(x,t)= \sum_{j=1}^\infty  (f,\varphi_j)\varphi_j\int_0^t (t-s)^{\alpha-1}E_{\alpha,\alpha}(-\lambda_j(t-s)^\alpha)\mu(s) \mathrm{d}s ,
\end{eqnarray}
where $E_{\alpha,\beta}(z)$ is the two-parameter Mittag-Leffler function: 
\begin{eqnarray*}
 E_{\alpha,\beta}(z):=\sum\limits_{k=0}^\infty\frac{z^k}{\Gamma(\alpha k+\beta)},\quad \forall z\in\mathds{C}.
\end{eqnarray*}

We define the source-to-solution operator $T: f\to u(x,T;f)$, where, in order to make clear the
dependence on $f$, we use $u(x,T;f)$ instead of $u(x,T)$ to denote the solution of
\eqref{dir-prob} associated to the source data $f$.
Then on the basis of \cite[Theorem 2.2]{Sak-Yam-2011}, the final state operator
$T:f\to u(x,T;f)$ is a continuous and linear mapping from $L^2(\Omega)$ into $L^2(\Omega)$ for any fixed $\mu\in L^\infty(0,T)$.
In general, $Tf=g$ doesn't necessarily have a solution for $g\in L^2(\Omega)$.
In this paper, we always assume that $g$ is attainable, i.e., $g\in\mathcal{R}(T)$,
where $\mathcal{R}(T)$ stands for the range of $T$.
The solution associated to an attainable $g$ is not necessarily unique except for
$\mu(t)\in C[0,T]$ satisfying $\mu(t)\geq\mu_0>0$ for all $t\geq 0$ \cite{Wei-Wang-2014}.
In fact, the operator equation $Tf=g$ has many solutions if the null space $\mathcal{N}(T)\neq\{0\}$.
Nevertheless, one might be interested in a specific solution satisfying additional requirements, such as the best
approximate solution \cite{Engl-Hanke-Beyvayer-1996}.
If $\mu\in L^\infty(0,T)$, we deduce from \eqref{dir-solution} and the identity \cite{Sak-Yam-2011}
\begin{eqnarray*}
 \frac{\mathrm{d}}{\mathrm{d}t}E_{\alpha,1}(-\lambda t^\alpha)
 =
 -\lambda t^{\alpha-1}E_{\alpha,\alpha}(-\lambda t^\alpha),\quad \lambda,\alpha,t>0
\end{eqnarray*}
that
\begin{align*}
(Tf)(x)=u(x,T;f)
&\leq \|\mu\|_{L^\infty(0,T)}\sum_{j=1}^\infty (f,\varphi_j)\varphi_j
\int_0^T (T-s)^{\alpha-1}E_{\alpha,\alpha}(-\lambda_j(T-s)^\alpha)\mathrm{d}s \\
&=\|\mu\|_{L^\infty(0,T)}\sum_{j=1}^\infty \frac{1-E_{\alpha,1}(-\lambda_jT^\alpha)}{\lambda_j} (f,\varphi_j)\varphi_j,
\end{align*}
It follows from the fact $(1-E_{\alpha,1}(-\lambda_jT^\alpha))/\lambda_j\to 0$ as $j\to\infty$
that reconstructing $f$ from $Tf=g$ is unstable.
This means that Problem (\textbf{IP}) is ill-posed.

Throughout the paper we will assume that $\mu(t)$ is an absolutely continuous function in $[0,T]$.
It is known \cite{Podlubny-1999} that
a such function
 possesses a (generalized) derivative in $L^1(0,T)$.

\subsection{$L^2$-TV regularization}

We briefly recall the space of functions having bounded total variation.
A scalar function $f\in L^1(\Omega)$ is said to be of bounded total variation if
\begin{eqnarray*}
TV(f):=\sup\limits_{\rho\in\mathcal{B}}\int_\Omega f\div\rho \mathrm{d}x<\infty \quad\text{for}\quad \mathcal{B}:= \left\{\rho~|~\rho\in C_c^1(\Omega)^d,~|\rho|_\infty\leq 1 \right\}
\end{eqnarray*}
where $|\rho|_\infty:=\sup\limits_{x\in \Omega}\left(\sum\limits_{i=1}^d|\rho_i(x)|^2\right)^{1/2}$,
$\div$ denotes the divergence operator, and $C_c^1(\Omega)$ is the space of
continuously differentiable $\mathds{R}^d$-valued functions with compact support in $\Omega$.
We use $BV(\Omega)$ to denote the space of all functions in $L^1(\Omega)$ of bounded total variation, i.e.,
\begin{eqnarray*}
BV(\Omega):=\{f\in L^1(\Omega)~|~TV(f)<\infty\},
\end{eqnarray*}
which is a Banach space equipped with the norm $\|f\|_{BV(\Omega)}:=\|f\|_{L^1(\Omega)} + TV(f)$
\cite{Bartels-2015,Giusti-1984,Ziemer-1989}.
Note that the space $BV(\Omega)$ is an extension of $W^{1,1}(\Omega)$ in the sense that $W^{1,1}(\Omega)\subset BV(\Omega)$ and $\|f\|_{BV(\Omega)}=\|f\|_{W^{1,1}(\Omega)}$ for all $f\in W^{1,1}(\Omega)$.

The following proposition is useful for our analysis, see \cite[Proposition 1.1]{Chavent-K-1997}.
\begin{prop}\label{prop-BVfun}
$(i)$~ $($Lower semicontinuity$)$ If $\{f_m\}_{m=1}^\infty\subset BV(\Omega)$ and $f_m\to f$ in $L^1(\Omega)$ as $m\to \infty$, then $f\in BV(\Omega)$ and
\begin{eqnarray*}
   &&  TV(f)\leq\liminf\limits_{m\to\infty}TV(f_m).
\end{eqnarray*}

$(ii)$~For every bounded sequence $\{f_m\}_{m=1}^\infty\subset BV(\Omega)$ there exists a subsequence $\{f_{m_k}\}_{k=1}^\infty$ and $f\in BV(\Omega)$ such that $f_{m_k}\to f$ in $L^p(\Omega)$ as $k\to\infty$, $p\in[1,\frac{d}{d-1})$, if $d\geq 2;$ and $f_{m_k}\to f$ in $L^p(\Omega)$ as $k\to\infty$, $p\in[1,\infty)$, if $d= 1$.

$(iii)$~$($Sobolev inequality$)$ For any $f\in BV(\Omega)$ and $1\leq p\leq d/(d-1)$, it holds
\begin{eqnarray}
   && \|f-\widetilde{f}\|_{L^p(\Omega)} \leq c TV(f), \quad\text{and}\quad \widetilde{f}:=\frac{1}{|\Omega|}\int_\Omega f\mathrm{d}x. \label{sob-inequ}
\end{eqnarray}
If $d=1$, then the $L^{d/(d-1)}$-norm is understood to be the $L^\infty$-norm.
\end{prop}

To deal with the ill-posedness of Problem (\textbf{IP}), we consider the following optimization problem
with both the $L^2$ and TV regularizations
\begin{align}
  \min\limits_{f\in BV(\Omega)\cap L^2(\Omega)} J(f),\ \
  \mbox{with }J(f):=\frac{1}{2}\|u(f)-g\|_{L^2(\Omega)}^2 + \frac{\beta}{2}\|f\|_{L^2(\Omega)}^2 + \gamma TV(f), 
  \label{min-L2TV-contin}
\end{align}
where
$\beta,\gamma>0$ are the regularization parameters.
For simplifying the notation we have used $u(f)$ to denote $u(x,T;f)$.

As pointed out in \cite{Chavent-K-1997,Tian-Yuan-2019}, it makes sense to add the $L^2$ quadratic regularization term in (\ref{min-L2TV-contin}).
First, the  $L^2$ term $\frac{\beta}{2}\|\cdot\|_{L^2(\Omega)}^2$ guarantees that
the functional $J^\delta(f)$
is strictly convex, so that the solution to \eqref{min-L2TV-contin} is unique.
Second, it provides a coercive term for the subspace of constant functions
which are in the kernel of the weak gradient operator
(and can be in the kernel of $T$).
Third, the space $BV(\Omega)\cap L^2(\Omega)$ endowed with the norm $\|f\|:=\|f\|_{L^2(\Omega)}+\|f\|_{BV(\Omega)}$ is a Banach space,
due to
the assumption that $\Omega$ is bounded
the fact that $BV(\Omega)\cap L^2(\Omega)$ and $BV(\Omega)$ are equivalent spaces if $d\leq 2$ (see \cite[Corollary 1.2]{Chavent-K-1997}).
Therefore, the optimization problem (\ref{min-L2TV-contin}) takes advantage of the favorable feature of the TV regularization.

For convenience, the regularized functional $J$ defined in \eqref{min-L2TV-contin} will be denoted by
$J^\delta$  in case $g$ is replaced by noisy data $g^\delta$, i.e.,
\be\label{Jdel}
J^\delta(f):=\frac{1}{2}\|u(f)-g^\delta\|_{L^2(\Omega)}^2 + \frac{\beta}{2}\|f\|_{L^2(\Omega)}^2 + \gamma TV(f).
\ee

The following existence and stability results can be directly proved by using
\cite[Theorems 2.1 and 3.1]{Chavent-K-1997}.
\begin{thm}\label{thm-conmin-exuniq}
For any $g\in L^2(\Omega)$, there exists a unique minimizer to the optimization problem \eqref{min-L2TV-contin}.
\end{thm}

\begin{thm}\label{thm-conmin-stable}
Let $\delta_i\to 0$ as $i\to\infty$. $\{g^{\delta_i}\}_{i=1}^\infty\subset L^2(\Omega)$ is a sequence of observations satisfying $g^{\delta_i}\rightarrow g$ in $L^2(\Omega)$ as $i\to\infty$. Let $\{f^{\delta_i}\}_{i=1}^\infty\subset BV(\Omega)\cap L^2(\Omega)$ denote the sequence of solutions to \eqref{min-L2TV-contin} corresponding to the observation
$g^{\delta_i}$. Then there exists a (not relabeled) subsequence of $\{f^{\delta_i}\}_{i=1}^\infty$ and $f^*\in BV(\Omega)\cap L^2(\Omega)$ such that
\begin{eqnarray*}
 \lim\limits_{i\to\infty}\|f^{\delta_i}-f^*\|_{L^2(\Omega)}=0,\quad\text{and}\quad  \lim\limits_{i\to\infty}TV(f^{\delta_i})=TV(f^*).
\end{eqnarray*}
Furthermore, $f^*$ is the unique solution to \eqref{min-L2TV-contin} corresponding to the noiseless observation $g$.
\end{thm}

\section{Some finite element approximation results for $BV$ space}\label{sec:x}

To establish our main results later, we first need to derive some approximation results for $BV$ space.

Let $\{\mathcal{T}_h\}_{0<h<1}$ be a family of regular triangulations of the domain $\Omega$, with the mesh parameter $h$ denoting the maximum diameter of the elements.
We introduce the standard finite element space $X_h^1$ of continuous piecewise linear functions:
\begin{eqnarray*}
  X_h^1=\left\{v_h\in C(\overline{\Omega})~|~ v_h|_K\in\mathds{P}_1(K),~\forall K\in\mathcal{T}_h \right\},
\end{eqnarray*}
where $\mathds{P}_1(K)$ denotes the space of linear polynomials on $K$.
Let $\Pi_h^1$ be the nodal interpolation operator \cite{Brenner-Scott-2007} from $C(\bar{\Omega})$ to $X_h^1$,
which allows the following interpolation error estimate:
\begin{eqnarray}
   && \|\Pi_h^1 v-v\|_{L^p(\Omega)} +h \|\nabla (\Pi_h^1 v-v)\|_{L^p(\Omega)}\leq c_\Pi h^2\|D^2v\|_{L^p(\Omega)},\quad \forall v\in W^{2,p}(\Omega), \label{est-nod-interp}
\end{eqnarray}
where $c_\Pi$ is a constant independent of $h$.
\begin{lem}\label{lem-stab-ipthm-1}
For any $f\in BV(\Omega)\cap L^p(\Omega), 1\leq p<\infty$,
there exists $\{\widehat{f}_h\}_{h>0}\subset X_h^1$ such that
\begin{align*}
\lim\limits_{h\to 0} \|\widehat{f}_{h}-f\|_{L^p(\Omega)}=0  \quad \text{and}\quad  \lim\limits_{h\to 0}TV(\widehat{f}_h)=TV(f).
\end{align*}
\end{lem}
\noindent \emph{Proof }
We want to prove that for any $f\in BV(\Omega)\cap L^p(\Omega)$ and $\delta>0$, there exists
$\{\widehat{f}_h\}_{h>0}\subset X_h^1$ such that
\be\label{te1}
\|\widehat{f}_h-f\|_{L^p(\Omega)}\le
\delta, \ \ \ TV(\widehat{f}_h)-TV(f)\leq \delta
\ee
for sufficiently small $h$.
To this end we first obtain from the density of the smooth functions in $BV(\Omega)\cap L^p(\Omega)$
\cite{Bartels-Nochetto-2014} that there exist $\varepsilon_0>0$ and $\{f_\varepsilon\}_{\varepsilon>0}\subset C^\infty(\Omega)\cap BV(\Omega)\cap L^p(\Omega)$ such that for all $\varepsilon\leq \varepsilon_0$, it holds
\be\label{app-BVp}
\ba{l}
\dps \|f_{\varepsilon}-f\|_{L^p(\Omega)}\leq {1\over 2}\delta, \qquad  \|\nabla f_{\varepsilon}\|_{L^1(\Omega)}\leq TV(f)+ {1\over 2}\delta, \\
 \|D^2f_{\varepsilon}\|_{L^1(\Omega)}\leq c\varepsilon^{-1}TV(f), \qquad \|D^2f_{\varepsilon}\|_{L^p(\Omega)}\leq c\varepsilon^{-2}\|f\|_{L^p(\Omega)},
\ea
\ee
where $c$ is a constant independent of $\delta$ and $\varepsilon$.
Let $\widehat{f}_h:=\Pi_h^1f_\varepsilon\in X_h^1$.
It follows from \eqref{est-nod-interp} and \eqref{app-BVp},
\begin{align*}
\|\widehat{f}_h-f\|_{L^p(\Omega)} &\leq \|\Pi_h^1f_\varepsilon-f_\varepsilon\|_{L^p(\Omega)} + \|f_\varepsilon-f\|_{L^p(\Omega)} \\
 & \leq c_\Pi h^2\|D^2f_\varepsilon\|_{L^p(\Omega)} + {1\over 2}\delta
 \leq c_\Pi ch^2\varepsilon^{-2}\|f\|_{L^p(\Omega)}+ {1\over 2}\delta,
\end{align*}
and
\begin{align*}
TV(\widehat{f}_h)-TV(f) &=\|\nabla \Pi_h^1f_{\varepsilon}\|_{L^1(\Omega)} - TV(f) \\
 & \leq \|\nabla (\Pi_h^1f_{\varepsilon}-f_\varepsilon)\|_{L^1(\Omega)}
 + \|\nabla f_{\varepsilon}\|_{L^1(\Omega)} - TV(f) \\
 & \leq c_\Pi c h\|D^2f_\varepsilon\|_{L^1(\Omega)} + {1\over 2}\delta
 \leq c_\Pi c h\varepsilon^{-1}TV(f) + {1\over 2}\delta.
\end{align*}
Then taking sufficiently small $h$ such that $c_\Pi ch^2\varepsilon^{-2}\|f\|_{L^p(\Omega)}\le {1\over 2}\delta$
and
$c_\Pi c h\varepsilon^{-1}TV(f)\le {1\over 2}\delta$,
we obtain \eqref{te1}.\\
A direct consequence of \eqref{te1} is $\|\widehat{f}_h-f\|_{L^p(\Omega)}\to 0$ and
$TV(\widehat{f}_h)\leq TV(f)$ as $h\to 0$.
Now we combine the latter with Proposition \ref{prop-BVfun}$(i)$ to get
\begin{eqnarray*}
   && TV(f)\leq \liminf\limits_{h\to 0}TV(\widehat{f}_h)\leq\limsup\limits_{h\to 0}TV(\widehat{f}_h)\leq TV(f),
\end{eqnarray*}
which implies $TV(\widehat{f}_h)\to TV(f)$ as $h\to 0$.
\hfill$\Box$

The following lemma provides a uniform approximation of the $BV$ space.

\begin{lem}\label{lem-stab-ipthm-2}
Assume that $f\in BV(\Omega)\cap L^\infty(\Omega)$. Then

$(i)$~for any $1\leq p<\infty$ and $h>0$, an element $\widehat{f}_{h}\in X_h^1$ exists such that
\begin{align}
\|\widehat{f}_{h}-f\|_{L^p(\Omega)}\leq ch^{\frac{1}{p+1}} \quad \text{and}\quad TV(\widehat{f}_h)-TV(f)\leq ch^{\frac{1}{p+1}}. \label{unapBV-inf-1}
\end{align}

$(ii)$~for any $1\leq p<\infty$ and $h>0$, an element $\widehat{f}_{h}\in X_h^1$ exists such that
\begin{align}
\|\widehat{f}_{h}-f\|_{L^p(\Omega)}\leq c(h|\ln h|)^{\frac{1}{p}} \quad \text{and}\quad \lim\limits_{h\to 0}TV(\widehat{f}_h)=TV(f). \label{unapBV-inf-1b}
\end{align}
\end{lem}

\noindent \emph{Proof }
It follows from \cite[Lemma 10.1]{Bartels-2015} that for $f\in BV(\Omega)$ and any $\varepsilon>0$,
there exists $f_{\varepsilon,h}\in X_h^1$ such that
\begin{align*}
\|f_{\varepsilon,h}-f\|_{L^1(\Omega)}\leq c(h^2\varepsilon^{-1}+\varepsilon) TV(f) \quad \text{and}\quad TV(f_{\varepsilon,h})\leq (1+ch\varepsilon^{-1}+c\varepsilon)TV(f).
\end{align*}
Moreover, it holds $\|f_{\varepsilon,h}\|_{L^\infty(\Omega)}\leq\|f\|_{L^\infty(\Omega)}$ under the assumption $f\in L^\infty(\Omega)$.
Let $\widehat{f}_h:=f_{\varepsilon,h}$, then
\begin{align*}
 \|\widehat{f}_h-f\|_{L^p(\Omega)}^p & =\int_{\Omega}|\widehat{f}_h-f|\cdot|\widehat{f}_h-f|^{p-1}\mathrm{d}x \nonumber\\
 & \leq \int_{\Omega}|\widehat{f}_h-f|\cdot(|\widehat{f}_h|+|f|)^{p-1}\mathrm{d}x \leq 2^{p-1}\|f\|_{L^\infty(\Omega)}^{p-1} \|\widehat{f}_h-f\|_{L^1(\Omega)}.
\end{align*}
Therefore,
\begin{align}
\|\widehat{f}_h-f\|_{L^p(\Omega)}\leq c(h^2\varepsilon^{-1}+\varepsilon)^{\frac{1}{p}} \quad \text{and}\quad TV(\widehat{f}_h)-TV(f)\leq c(h\varepsilon^{-1}+\varepsilon),
\label{lem-stab-ipthm-2-eq1}
\end{align}
where $c$ may depend on $f$. \\
Obviously, we obtain \eqref{unapBV-inf-1} by taking $\varepsilon\sim h^{\frac{p}{p+1}}$ in \eqref{lem-stab-ipthm-2-eq1}.
Then we take $\varepsilon\sim h|\ln h|$ in \eqref{lem-stab-ipthm-2-eq1} to get the first estimate of \eqref{unapBV-inf-1b}
and
\begin{eqnarray*}
TV(\widehat{f}_h)-TV(f)\leq c(\ln h)^{-1}.
\end{eqnarray*}
This last inequality allows to derive the second estimate of \eqref{unapBV-inf-1b} by following the same lines as
Lemma \ref{lem-stab-ipthm-1}.
\hfill$\Box$

\begin{rem}
Optimizing the convergence rates of the estimates in \eqref{lem-stab-ipthm-2-eq1} simultaneously leads to the choice $\varepsilon\sim h^{p/(p+1)}$ as $h\to 0$.
However, in the subsequent analysis, we are going to see that the convergence result of the proposed algorithm can be
guaranteed by the simple convergence of $TV(\widehat{f}_h)$; see Theorems \ref{thm-convfm-ipthm} and \ref{thm-stab-data}. Therefore, as suggested by \cite{Hinze-Kalt-2018} we will take
$\varepsilon\sim h|\ln h|$, which leads to the convergence rate
$O(h^r)$ as $h\to 0$ for $r\in(0,1)$.
\end{rem}

For $1\leq p<\infty, 0<s\leq 1$, we introduce the Lipschitz space $\mathrm{Lip}(s,L^p(\Omega))$
\cite[Chapter 2]{DeVore-Lorentz-1993}, which consists of all functions $v\in L^p(\Omega)$ such that
\begin{eqnarray*}
   && |v|_{\mathrm{Lip}(s,L^p(\Omega))} :=\sup\limits_{t>0}\big\{ t^{-s}\omega(v,t)_p \big\}<\infty,\quad\text{with}\quad \omega(v,t)_p=\sup\limits_{|y|\leq t}\big( \int_\Omega|v(x+y)-v(x)|^p \mathrm{d}x \big)^{1/p}.
\end{eqnarray*}
We note that $BV(\Omega)\subset \mathrm{Lip}(1,L^1(\Omega))$. Furthermore, it is known; see
\cite[Lemma 1.1 and Lemma 1.6]{Wang-Lucier-2011} and \cite[Section 2.4]{Bartels-2012}, that
for $f\in\mathrm{Lip}(s,L^p(\Omega))$ and any $\varepsilon>0$, there exist
$\{f_\varepsilon\}_{\varepsilon>0}\subset C^\infty(\Omega)\cap BV(\Omega)\cap L^p(\Omega)$ such that
\begin{eqnarray}
   && \|f_{\varepsilon}-f\|_{L^p(\Omega)}\leq c\varepsilon^s, \qquad  \|\nabla f_{\varepsilon}\|_{L^1(\Omega)}\leq TV(f), \qquad \|D^2f_{\varepsilon}\|_{L^1(\Omega)}\leq c\varepsilon^{-1}TV(f) \label{app-BV-Lip1},
\end{eqnarray}
and
\begin{eqnarray}
   && \omega(f_\varepsilon,t)_p\leq c\omega(f,t)_p,\ \ \ \forall t>0.
\label{app-BV-Lip2}
\end{eqnarray}

\begin{lem}\label{lem-stab-ipthm-X}
Assume that $f\in BV(\Omega)\cap \mathrm{Lip}(s,L^p(\Omega))$ for some $1\leq p<\infty$ and $0<s\leq 1$. Then
there exits $\widehat{f}_{h}\in X_h^1$ such that
\begin{align}
\|\widehat{f}_{h}-f\|_{L^p(\Omega)}\leq ch^{\frac{s}{s+1}} \quad \text{and}\quad TV(\widehat{f}_h)-TV(f)\leq ch^{\frac{s}{s+1}}
\label{unapBV-Lip-1}
\end{align}
or
\begin{align}
\|\widehat{f}_{h}-f\|_{L^p(\Omega)}\leq c(h|\ln h|)^{s} \quad \text{and}\quad \lim\limits_{h\to 0}TV(\widehat{f}_h)=TV(f). \label{unapBV-Lip-2}
\end{align}
\end{lem}

\noindent \emph{Proof } For $f\in BV(\Omega)\cap \mathrm{Lip}(s,L^p(\Omega)), 1\leq p<\infty, 0<s\leq 1$,
let $\{f_\varepsilon\}_{\varepsilon>0}$ be a function in $C^\infty(\Omega)\cap BV(\Omega)\cap L^p(\Omega)$ satisfying
\eqref{app-BV-Lip1} and \eqref{app-BV-Lip2}.
Define $\widehat{f}_h:={}^{q}\Pi_h^1f_\varepsilon\in X_h^1$, where
${}^{q}\Pi_h^1$ is the quasi-interpolation operator; see, e.g., \cite[Definition 3.8]{Bartels-2015} for the definition.
Then we use an approximation result established in \cite[Theorem 7.3]{DeVore-Lorentz-1993}, i.e.,
\begin{eqnarray*}
   && \| {}^{q}\Pi_h^1f_\varepsilon -f_\varepsilon\|_{L^p(\Omega)}\leq c \omega(f_\varepsilon,h)_p
\end{eqnarray*}
to obtain
\begin{align*}
\|\widehat{f}_h-f\|_{L^p(\Omega)}
& \leq \| {}^{q}\Pi_h^1f_\varepsilon -f_\varepsilon \|_{L^p(\Omega)} + \|f_\varepsilon-f \|_{L^p(\Omega)} \\
& \leq c\omega(f_\varepsilon,h)_p + c\varepsilon^s  \leq c\omega(f,h)_p +c\varepsilon^s \leq c(h^s+\varepsilon^s).
\end{align*}
On the other hand, we have \cite[Lemma 10.1]{Bartels-2015}
\begin{eqnarray*}
   && TV(\widehat{f_h})=\|\nabla {}^{q}\Pi_h^1f_\varepsilon\|_{L^1(\Omega)}\leq(1+ch\varepsilon^{-1}+c\varepsilon)TV(f).
\end{eqnarray*}
Hence
\begin{align*}
TV(\widehat{f_h})-TV(f)\leq c(h\varepsilon^{-1}+\varepsilon).
\end{align*}
The rest of the analysis is the same as the proof of Lemma \ref{lem-stab-ipthm-2}.
\hfill$\Box$

\begin{rem}
$(i)$~Since $BV(\Omega)\subset \mathrm{Lip}(1,L^1(\Omega))$, Lemma \ref{lem-stab-ipthm-X} implies that for $f\in BV(\Omega)$,
an element $\widehat{f}_{h}\in X_h^1$ exists such that
\begin{align*}
\|\widehat{f}_{h}-f\|_{L^1(\Omega)}\leq ch|\ln h| \quad \text{and}\quad \lim\limits_{h\to 0}TV(\widehat{f}_h)=TV(f),
\end{align*}
which can be found in \cite[Lemma 4.6]{Hinze-Kalt-2018}.

$(ii)$~By the Sobolev inequality \eqref{sob-inequ} and \cite{Bartels-2012,Wang-Lucier-2011}, we have $BV(\Omega)\cap L^\infty(\Omega)=BV(\Omega)\subset \mathrm{Lip}(1,L^1(\Omega))$ if $d=1$, and $BV(\Omega)\cap L^\infty(\Omega)\subset \mathrm{Lip}(1/2,L^2(\Omega))$ if $d=2$.
In these cases, Lemma \ref{lem-stab-ipthm-2} and Lemma \ref{lem-stab-ipthm-X} give the same result.
\end{rem}

\section{Discretization of the direct problem and convergence}\label{sec:3}

This section is devoted to propose and analyze a full discretization method for the direct problem (\ref{dir-prob}).
The approximation is based on a finite element discretization in space and a finite
difference scheme in time.

\subsection{Discretization of the direct problem}

We begin with a space semi-discrete scheme for solving (\ref{dir-prob}) by finite element method.
Let $V_h:=H_0^1\cap X_h^1$.
The semi-discrete problem to (\ref{dir-prob}) reads: find $u_h(t)\in V_h$, such that
\begin{eqnarray}
 u_h(0)=0,\quad (\partial_t^\alpha u_h(t), v_h ) + (\nabla u_h(t),\nabla v_h) = \mu(t)(f,v_h),\quad \forall v_h\in V_h,\quad 0<t\leq T. \label{semi-dis-dirp}
\end{eqnarray}
A straightforward application of \cite[Theorem 3.7]{Jin-Lazarov-2015} gives the following error estimate.
\begin{lem}\label{lem-forp-semi-erro}
Let $f\in L^2(\Omega)$, and $u$ and $u_h$ be the solution of \eqref{dir-prob} and \eqref{semi-dis-dirp}, respectively. Then there holds
\begin{eqnarray*}
 && \|u_h(t)-u(t)\|_{L^2(\Omega)}+h\|\nabla (u_h(t)-u(t))\|_{L^2(\Omega)} \leq ch^2|\ln h|^2 \|\mu\|_{L^\infty(0,T)} \|f\|_{L^2(\Omega)},\quad \forall t>0.
\end{eqnarray*}
\end{lem}

Then we turn to describe the time discretization, which makes use of the popular L1 scheme to discretize the time-fractional derivative
 \cite{Sun-Wu-2006,lin-xu-2007}.
Let $t_k:=k\tau$, $k=0,1,\cdots,K_\tau$, where $\tau:=T/K_\tau$ is the time step.
The fractional derivative \eqref{det-caputo-1} at $t=t_{k+1}$ is approximated by
\begin{align}
\partial_t^\alpha u(x,t_{k+1}) &= \frac{1}{\Gamma(1-\alpha)} \sum\limits_{j=0}^k \int_{t_j}^{t_{j+1}} (t_{k+1}-s)^{-\alpha}\frac{\partial u(x,s)}{\partial s} \mathrm{d}s  \nonumber\\
& =\frac{1}{\Gamma(1-\alpha)} \sum\limits_{j=0}^k \frac{u(x,t_{j+1})-u(x,t_j)}{\tau} \int_{t_j}^{t_{j+1}} (t_{k+1}-s)^{-\alpha} \mathrm{d}s + r_\tau^{k+1} \nonumber\\
& =\frac{1}{\Gamma(2-\alpha)\tau^\alpha} \sum\limits_{j=0}^k b_j(u(x,t_{k-j+1})-u(x,t_{k-j})) + r_\tau^{k+1} \nonumber\\
& =\frac{1}{\Gamma(2-\alpha)\tau^\alpha} \sum\limits_{j=0}^{k+1}\zeta_j u(x,t_{k-j+1}) + r_\tau^{k+1} \nonumber\\
& :=L_{\tau}^\alpha u(x,t_{k+1}) + r_\tau^{k+1}, \label{L1-scheme}
\end{align}
where $b_j:=(j+1)^{1-\alpha}-j^{1-\alpha}$ for $j=0,1,\cdots,k$, $r_\tau^{k+1}$ is the local truncation error, and the coefficients $\zeta_j$ in the discrete operator $L_{\tau}^\alpha$ are given by
\begin{eqnarray*}
   && \zeta_j=\left\{\begin{array}{ll}
                       1, & j=0, \\
                       (j+1)^{1-\alpha}-2j^{1-\alpha}+(j-1)^{1-\alpha}, & j=1,\cdots,k, \\
                       k^{1-\alpha}-(k+1)^{1-\alpha}, & j=k+1.
                     \end{array}\right.
\end{eqnarray*}
Note that $b_j$ and $\zeta_j$ satisfy
\begin{eqnarray}
    && 1=b_0>b_1>\cdots >b_k>0,\qquad b_k\to 0\ \text{as}\  k\to\infty, \qquad \sum\limits_{j=0}^{k-1}(b_j-b_{j+1}) + b_k =1.\label{rel-zetaj} \\
    && \zeta_0>0>\zeta_{k+1}>\zeta_k>\cdots>\zeta_1,\qquad \sum\limits_{j=0}^{n}\zeta_j\geq 0\ \ \text{for~any}~n=1, \cdots, k+1.
    \nonumber
\end{eqnarray}
We consider the fully discrete problem: find $u_h^{k+1}\in V_h$, such that
\begin{eqnarray}
   u_h^0=0,\quad (L_{\tau}^\alpha u_h^{k+1},v_h) + (\nabla u_h^{k+1},\nabla v_h) = \mu^{k+1}(f,v_h),\quad\forall v_h\in V_h,
   \ \forall k\geq 0, \label{fuly-dis-dirp-L}
\end{eqnarray}
where $u_h^k$ is an approximation to $u_h(t_k)$. Equivalently, \eqref{fuly-dis-dirp-L} can be rewritten as
\begin{eqnarray}
  u_h^0=0,\quad (u_h^{k+1},v_h) + \eta(\nabla u_h^{k+1},\nabla v_h) = \sum\limits_{j=0}^{k-1}(b_j-b_{j+1}) (u_h^{k-j},v_h) + \eta \mu^{k+1}(f,v_h),\quad \forall v_h\in V_h, \ \forall k\geq 0,
\label{fuly-dis-dirp}
\end{eqnarray}
where $\eta:=\Gamma(2-\alpha)\tau^\alpha$.
It was shown in \cite{Sun-Wu-2006,lin-xu-2007} that the time discretization error is of order ${2-\alpha}$,
provided that the exact solution
is twice continuously differentiable in time. However, in view of the smoothing property of the subdiffusion equation,
this regularity condition is restrictive, which
does not hold even for the homogeneous problem with a smooth initial data.
\cite{Jin-IMA-2016,Jin-Li-2018} revisited the error
analysis of the L1 scheme, and established an $O(\tau)$ convergence rate for solutions of low regularity, which can be directly
applied to our problem to give the following error estimate under the condition
\be\label{te2}
 \int_0^{t_k}(t_k-s)^{\alpha-1}\|\mu'(s)f\|_{L^2(\Omega)}\mathrm{d}s<\infty
\ee

\begin{lem}\label{lem-err-time}
If $u_h(t)$ and $u_h^k$ are the solutions of \eqref{semi-dis-dirp} and \eqref{fuly-dis-dirp}, respectively. Then there holds
\begin{eqnarray*}
   && \|u_h^k-u_h(t_k)\|_{L^2(\Omega)} \leq ct_k^{\alpha-1}\tau \|\mu\|_{L^\infty(0,T)} \|f\|_{L^2(\Omega)},\qquad \forall k\geq 1.
\end{eqnarray*}
\end{lem}
Notice that \eqref{te2} holds under the absolute continuity assumption on $\mu(t)$ since
\begin{align*}
   \int_0^{t_k}(t_k-s)^{\alpha-1}\|\mu'(s)f\|_{L^2(\Omega)}\mathrm{d}s &\leq \|f\|_{L^2(\Omega)}\int_0^{t_k}(t_k-s)^{\alpha-1}|\mu'(s)|\mathrm{d}s\\
   & \leq t_k^{\alpha-1}\|f\|_{L^2(\Omega)}\int_0^{t_k}|\mu'(s)|\mathrm{d}s \leq ct_k^{\alpha-1}\|\mu'\|_{L^1(0,T)}\|f\|_{L^2(\Omega)}<\infty.
\end{align*}

The properties of the finite difference operator $L_\tau$ plays a key role in the stability and convergence analysis
of the discrete solution $u_h^k$.
It can be found in \cite[Lemma 3.1]{Jiang-Liu-2020} and \cite[Lemma 3.1]{Wang-Xiao-Zou-2020}
that the operator $L_\tau^\alpha$ satisfies: for any function sequence $v_h^k\in V_h, k\geq 1$,
\be\label{La}
L_\tau^\alpha \|v_h^k\|_{L^2(\Omega)}^2 \leq 2(v_h^k, L_\tau^\alpha v_h^k),\qquad \forall k\geq 1.
\ee
Furthermore, a simplified and improved result of \cite[Lemma 3.1]{Li-Liao-Sun-2018} shows:
if $ L_\tau^\alpha \omega^k \leq \theta^k$ for all $k\geq 1$, then
\be\label{La2}
   && \omega^k \leq 2\left(\omega^0 + \frac{t_k^\alpha}{\Gamma(1+\alpha)}\max\limits_{0\leq j\leq k}\theta^j \right),\qquad \forall k\geq 1.
\ee
\begin{thm}\label{thm-uncstable-fuly}
The fully discrete problem \eqref{fuly-dis-dirp} is stable with respect to the source function $f$ in the sense that for all $h$ and
$\tau>0$, there exists a constant $c$ independent of $\tau$ and $h$, such that
\begin{eqnarray}
\|u_h^k\|_{H^1(\Omega)}\leq c \|f\|_{L^2(\Omega)}, \qquad \forall k\geq 1. \label{disA-bound}
\end{eqnarray}
\end{thm}
\noindent\emph{Proof} First we will prove by induction that
\begin{eqnarray}
 && \|u_h^j\|_{L^2(\Omega)} \leq \frac{\eta}{b_{j-1}}\|\mu\|_{L^\infty(0,T)}\|f\|_{L^2(\Omega)}, \qquad \forall j\geq 1. \label{thm-full-stab-eq1}
\end{eqnarray}
When $j=1$, we have
\begin{eqnarray*}
   && (u_h^1,v_h) + \eta(\nabla u_h^1,\nabla v_h)=\eta \mu^1(f,v_h),\qquad \forall v_h\in V_h.
\end{eqnarray*}
Taking $v_h=u_h^1$ and using Cauchy-Schwarz inequality, we obtain immediately
\begin{eqnarray*}
   && \|u_h^1\|_{L^2(\Omega)} \leq \eta|\mu^1|\|f\|_{L^2(\Omega)} \leq \eta \|\mu\|_{L^\infty(0,T)}\|f\|_{L^2(\Omega)}=\frac{\eta}{b_0}\|\mu\|_{L^\infty(0,T)}\|f\|_{L^2(\Omega)}.
\end{eqnarray*}
Assuming \eqref{thm-full-stab-eq1} holds for all $j\leq k$, we want to show \eqref{thm-full-stab-eq1} for $j=k+1$. Taking $v_h=u_h^{k+1}$ in (\ref{fuly-dis-dirp}) gives
\begin{align*}
\|u_h^{k+1}\|_{L^2(\Omega)}^2 + \eta\|\nabla u_h^{k+1}\|_{L^2(\Omega)}^2 &= \sum\limits_{j=0}^{k-1}(b_j-b_{j+1}) (u_h^{k-j},u_h^{k+1}) + \eta \mu^{k+1}(f,u_h^{k+1})\\
& \leq  \sum\limits_{j=0}^{k-1}\big[(b_j-b_{j+1}) \|u_h^{k-j}\|_{L^2(\Omega)}+\eta\|\mu\|_{L^\infty(0,T)}\|f\|_{L^2(\Omega)}\big]\|u_h^{k+1}\|_{L^2(\Omega)}.
\end{align*}
Hence, by using the induction assumption and then (\ref{rel-zetaj}), we have
\begin{align*}
\|u_h^{k+1}\|_{L^2(\Omega)} &\leq \sum\limits_{j=0}^{k-1}(b_j-b_{j+1}) \|u_h^{k-j}\|_{L^2(\Omega)} + \eta \|\mu\|_{L^\infty(0,T)}\|f\|_{L^2(\Omega)}\\
& \leq  \sum\limits_{j=0}^{k-1}\frac{b_j-b_{j+1}}{b_{k-j-1}}\eta \|\mu\|_{L^\infty(0,T)}\|f\|_{L^2(\Omega)} + \eta \|\mu\|_{L^\infty(0,T)}\|f\|_{L^2(\Omega)}\\
& \leq  \sum\limits_{j=0}^{k-1}\frac{b_j-b_{j+1}}{b_{k}}\eta \|\mu\|_{L^\infty(0,T)}\|f\|_{L^2(\Omega)} + \eta \|\mu\|_{L^\infty(0,T)}\|f\|_{L^2(\Omega)} \\
& =\frac{\eta}{b_k}\|\mu\|_{L^\infty(0,T)}\|f\|_{L^2(\Omega)}.
\end{align*}
This proves (\ref{thm-full-stab-eq1}).
\\
On the other hand, noticing $k^{-\alpha}b_{k-1}^{-1}\leq\frac{1}{1-\alpha},\ \forall k\geq 1$ \cite{lin-xu-2007},
we have
\begin{eqnarray*}
   \frac{\eta}{b_{k-1}}=\Gamma(2-\alpha)T^\alpha K_\tau^{-\alpha}b_{k-1}^{-1}\leq\Gamma(2-\alpha)T^\alpha K_\tau^{-\alpha}b_{K_\tau-1}^{-1}
   \leq\frac{\Gamma(2-\alpha)T^\alpha}{1-\alpha},\qquad \forall k\geq 1.
\end{eqnarray*}
Consequently we obtain
\begin{eqnarray*}
   && \|u_h^k\|_{L^2(\Omega)}\leq c \|f\|_{L^2(\Omega)},\qquad \forall k\geq 1,
\end{eqnarray*}
where $c>0$ is constant independent of $\tau$ and $h$.
\\
Now we turn to prove the stability in $H^1$ semi-norm. Taking $v_h=2L_\tau^\alpha u_h^k$ in \eqref{fuly-dis-dirp-L} yields
\begin{eqnarray}
   &&  2(L_{\tau}^\alpha u_h^{k},L_{\tau}^\alpha u_h^{k}) + 2(\nabla u_h^{k},\nabla L_{\tau}^\alpha u_h^{k}) = 2\mu^k(f,L_{\tau}^\alpha u_h^{k}),\qquad \forall k\geq 1. \label{thm-full-stab-eq2}
\end{eqnarray}
By applying \eqref{La} to the second term in the left hand side, we obtain
\begin{eqnarray*}
   &&  2\|L_{\tau}^\alpha u_h^{k}\|_{L^2(\Omega)}^2 + L_\tau^\alpha\|\nabla u_h^k\|_{L^2(\Omega)}^2
    \leq \|\mu\|_{L^\infty(0,T)}^2\|f\|_{L^2(\Omega)}^2 + \|L_{\tau}^\alpha u_h^{k}\|_{L^2(\Omega)}^2 ,\qquad\forall k\geq 1.
\end{eqnarray*}
Hence we obtain
\begin{eqnarray*}
   && L_\tau^\alpha\|\nabla u_h^k\|_{L^2(\Omega)}^2  \leq \|\mu\|_{L^\infty(0,T)}^2\|f\|_{L^2(\Omega)}^2,\qquad\forall k\geq 1.
\end{eqnarray*}
Then, it follows from \eqref{La2} that
\begin{align*}
\|\nabla u_h^k\|_{L^2(\Omega)}^2 &\leq 2\left( \|\nabla u_h^0\|_{L^2(\Omega)}^2 + \frac{t_k^\alpha}{\Gamma(1-\alpha)} \|\mu\|_{L^\infty(0,T)}^2\|f\|_{L^2(\Omega)}^2 \right) \leq c\|f\|_{L^2(\Omega)}^2.
\end{align*}
 This completes the proof.
\hfill$\Box$

\begin{thm}\label{thm-fuly-erest}
Let $u$ and $u_h^k$ be the solutions of \eqref{dir-prob} and \eqref{fuly-dis-dirp}, respectively. Then
\begin{eqnarray*}
\|u(t_k)-u_h^k\|_{L^2(\Omega)} \leq c\|\mu\|_{L^\infty(0,T)}\|f\|_{L^2(\Omega)}(t_k^{\alpha-1}\tau  + h^2|\ln h|^2),\qquad \forall k\geq 1.
\end{eqnarray*}
\end{thm}
\noindent\emph{Proof} It follows directly from Lemma \ref{lem-forp-semi-erro} and Lemma \ref{lem-err-time} and the triangle inequality.
\hfill$\Box$

\subsection{Finite element approximation of the inverse problem}

Motivated by the fact that the piecewise constant and piecewise affine globally continuous finite element spaces
are dense in $BV(\Omega)$ with respect to weak$^*$ convergence in $BV(\Omega)$ \cite{Bartels-2012}, and
$BV(\Omega)$ is continuously embedded into $L^r(\Omega)$ with $1\leq r\leq d/(d-1)$ \cite[Theorem 10.1.4]{Attouch-B-M-2006},
we propose to use the finite element space $X_h^1$ to approximate the source function $f$ in (\ref{min-L2TV-contin}).
Precisely, we consider the following finite element discrete version of the optimization problem (\ref{min-L2TV-contin}):
\begin{align}
 \min\limits_{f\in X_h^1} J_{h,\tau}(f),\qquad J_{h,\tau}(f):=\frac{1}{2}\|u_h^{K_{\tau}}(f)-g\|_{L^2(\Omega)}^2 + \frac{\beta}{2}\|f\|_{L^2(\Omega)}^2 + \gamma TV(f),
 \label{min-L2TV-disc-0}
\end{align}
where $u_h^{K_{\tau}}(f)$ is the solution of (\ref{fuly-dis-dirp}) corresponding to the
source term $f$.

\begin{thm}\label{thm-discmin-exuniq}
Let $\beta,\gamma>0$ be the fixed regularization parameters.
For every fixed $h,\tau>0$  and $g\in L^2(\Omega)$,
there exists a unique minimizer to the discrete optimization problem \eqref{min-L2TV-disc-0}.
\end{thm}

\noindent \emph{Proof }First, the non-negativity of $J_{h,\tau}(f)$ means that the infimum of $J_{h,\tau}(f)$ is finite over $X_h^1$.
Thus there exists a sequence $\{f_m\}_{m=1}^\infty\subset X_h^1$ such that
\begin{eqnarray}
\lim\limits_{m\to\infty}J_{h,\tau}(f_m)=\inf\limits_{f\in X_h^1}J_{h,\tau}(f).  \label{thm-disex-eq1}
\end{eqnarray}
Moreover, it follows from the definition of $J_{h,\tau}$ that $\{f_m\}_{m=1}^\infty$ is uniformly bounded in $L^2(\Omega)\cap BV(\Omega)$.
By Proposition \ref{prop-BVfun} $(ii)$, there exists a (not relabeled) subsequence of $\{f_m\}_{m=1}^\infty$ and an element $f_{h,\tau}\in X_h^1$ such that
\begin{eqnarray*}
   f_m\rightharpoonup f_{h,\tau}\ \text{in}\ L^2(\Omega), \quad f_m\to f_{h,\tau}\ \text{in}\ L^1(\Omega),  \quad \text{as }m\to\infty,
\end{eqnarray*}
and
\begin{eqnarray}
  TV(f_{h,\tau})\leq\liminf\limits_{m\to\infty}TV(f_m). \label{thm-disex-eq2}
\end{eqnarray}
Next we want to show that $f_{h,\tau}$ is a minimizer of (\ref{min-L2TV-disc-0}).
It follows from Theorem \ref{thm-uncstable-fuly} that for each $k=1,2,\cdots,K_\tau$, the sequence $\{u_h^k(f_m)\}_{m=1}^\infty$ is uniformly bounded in $H_0^1(\Omega)$
with respect to $m$.
Since all norms are equivalent in finite-dimensional space, and weak convergence means
strong convergence.
As a result there exists a (not relabeled) subsequence of
$\{u_h^k(f_m)\}_{m=1}^\infty$ and an element $\widehat{u}_h^k\in V_h$ such that
\begin{eqnarray}
 u_h^k(f_m)\to \widehat{u}_h^k\ \text{in}\ H^1(\Omega),\quad f_m\to f_{h,\tau}\ \text{in}\ L^2(\Omega),\quad  \text{as } m\to\infty.
\label{thm-disex-eq4}
\end{eqnarray}
We can claim $\widehat{u}_h^k=u_h^k(f_{h,\tau})$. In fact, by (\ref{fuly-dis-dirp}), $\{u_h^k(f_m)\}_{k=1}^{K_\tau}$ satisfies $u_h^0(f_m)=0$ and
\begin{eqnarray*}
   (u_h^{k+1}(f_m),v_h) + \eta(\nabla u_h^{k+1}(f_m),\nabla v_h) = \sum\limits_{j=0}^{k-1}(b_j-b_{j+1}) (u_h^{k-j}(f_m),v_h) + \eta \mu^{k+1}(f_m,v_h),\quad \forall v_h\in V_h
\end{eqnarray*}
for all $m$. Passing to the limit as $m\to\infty$ in the above equation gives immediately
$\widehat{u}_h^k=u_h^k(f_{h,\tau})$.
\\
Combining (\ref{thm-disex-eq1})-(\ref{thm-disex-eq4}) and the lower semi-continuity of the $L^2$-norm, we can get
\begin{align*}
J_{h,\tau}(f_{h,\tau}) &= \frac{1}{2}\|u_h^{K_{\tau}}(f_{h,\tau})-g\|_{L^2(\Omega)}^2 + \frac{\beta}{2}\|f_{h,\tau}\|_{L^2(\Omega)}^2 + \gamma TV(f_{h,\tau}) \\
                &= \frac{1}{2}\|\widehat{u}_h^{K_{\tau}}-g\|_{L^2(\Omega)}^2 + \frac{\beta}{2}\|f_{h,\tau}\|_{L^2(\Omega)}^2 + \gamma TV(f_{h,\tau}) \\
                &\leq\liminf\limits_{m\to\infty}\left\{ \frac{1}{2}\|u_h^{K_{\tau}}(f_m)-g\|_{L^2(\Omega)}^2 + \frac{\beta}{2}\|f_m\|_{L^2(\Omega)}^2 + \gamma TV(f_m) \right\}\\
                & =\liminf\limits_{m\to\infty} J_{h,\tau}(f_m)=\inf\limits_{f\in X_h^1}J_{h,\tau}(f).
\end{align*}
Therefore $f_{h,\tau}$ is indeed a minimizer to the discrete optimization problem (\ref{min-L2TV-disc-0}).
\\
Finally, it is easy to see that $u_h^{K_\tau}(f)$ is linear with respect to $f$.
The uniqueness of $f_{h,\tau}$ is a direct consequence of the strict convexity of $J_{h,\tau}(f)$ with respect to $f$.
\hfill$\Box$

To prove the convergence of finite element approximation for \eqref{min-L2TV-disc-0},
we need the following strong convergence lemma.
\begin{lem}\label{lem-stab-ipthm-3}
Let $\{h_m\}_{m=1}^\infty$ and $\{\tau_m\}_{m=1}^\infty$ be any positive sequences such that $h_m\to 0$ and $\tau_m\to 0$ as $m\to \infty$.
Let the sequence $\{f_m\}_{m=1}^\infty\subset BV(\Omega)\cap L^2(\Omega)$ converges weakly to $f\in BV(\Omega)\cap L^2(\Omega)$ in $L^2$-norm as $m\to \infty$. Then for any $g\in L^2(\Omega)$, it holds
\begin{eqnarray*}
   && \lim\limits_{m\to \infty} \|u_{h_m}^{K_{\tau_m}}(f_{m})-g\|_{L^2(\Omega)}=\|u(f)-g\|_{L^2(\Omega)}.
\end{eqnarray*}
\end{lem}
\noindent \emph{Proof }
Note that $f_m\rightharpoonup f$ implies $\{\|f_m\|_{L^2(\Omega)}\}_{m=1}^\infty$ is bounded. By using the a priori estimates of $u_{h_m}^{K_{\tau_m}}(f_m)$ and $u(f)$, we have
\begin{align*}
\left| \|u_{h_m}^{K_{\tau_m}}(f_{m})-g\|_{L^2(\Omega)}-\|u(f)-g\|_{L^2(\Omega)} \right|
& = \left| \big(u_{h_m}^{K_{\tau_m}}(f_{m})-u(f),u_{h_m}^{K_{\tau_m}}(f_{m})+u(f)-2g\big) \right| \\
&\leq \|u_{h_m}^{K_{\tau_m}}(f_{m})-u(f)\|_{L^2(\Omega)}
\|u_{h_m}^{K_{\tau_m}}(f_{m})+u(f)-2g\|_{L^2(\Omega)}\\
&\leq c \|u_{h_m}^{K_{\tau_m}}(f_{m})-u(f)\|_{L^2(\Omega)}.
\end{align*}
Then it suffices to show that $\|u_{h_m}^{K_{\tau_m}}(f_{m})-u(f)\|_{L^2(\Omega)}\to 0$ as $m\to\infty$.
We split the error $u_{h_m}^{K_{\tau_m}}(f_{m})-u(f)$ as
\begin{eqnarray*}
   && u_{h_m}^{K_{\tau_m}}(f_{m})-u(f)=u_{h_m}^{K_{\tau_m}}(f_{m})-u_{h_m}^{K_{\tau_m}}(f)
    + u_{h_m}^{K_{\tau_m}}(f)-u(f).
\end{eqnarray*}
First, it follows from Theorem \ref{thm-fuly-erest} that the second part, i.e.,
$\|u_{h_m}^{K_{\tau_m}}(f)-u(f)\|_{L^2(\Omega)}\to 0$ as $m\to \infty$.
\\
Now, for any $m>0$, let
$\xi_{h_m}^k:=u_{h_m}^{k}(f_{m})-u_{h_m}^{k}(f)$ for $k=0,1,\cdots,K_{\tau_m}$.
We deduce from \eqref{fuly-dis-dirp} that $\xi_{h_m}^k$ satisfies $\xi_{h_m}^0=0$, and
for all $k\geq 1$,
\begin{eqnarray}
 (\xi_{h_m}^{k+1},v_{h_m}) + \eta(\nabla \xi_{h_m}^{k+1},\nabla v_{h_m}) = \sum\limits_{j=0}^{k-1}(b_j-b_{j+1}) (\xi_{h_m}^{k-j},v_{h_m}) + \eta \mu^{k+1}(f_m-f,v_{h_m}),\quad \forall v_{h_m}\in V_{h_m}. \label{lem-stab-ipthm-3-eq1}
\end{eqnarray}
According to Theorem \ref{thm-uncstable-fuly} and the weak convergence of $f_m$ to $f$
in $L^2(\Omega)$, we have
\begin{eqnarray*}
 \|\xi_{h_m}^k\|_{H^1(\Omega)}
 \leq c \|f_m-f\|_{L^2(\Omega)} \leq c,\ \ \text{ for } m \text{ large enough}.
\end{eqnarray*}
Hence there exist a (not relabeled) subsequence of $\{\xi_{h_m}^k\}_{m=1}^\infty$ and $\xi^k\in H^1(\Omega)$ such that
\begin{eqnarray}
  \xi_{h_m}^k\rightharpoonup \xi^k\ \text{in}\ H^1(\Omega), \quad
 \xi_{h_m}^k\to \xi^k\ \text{in}\ L^2(\Omega),\quad\text{as } m\to \infty,
 \label{lem-stab-ipthm-3-eq2}
\end{eqnarray}
for all $k\geq 1$.
\\
Let
$R_h:H_0^1(\Omega)\to V_h$ is the Ritz projection defined by: $\forall \psi\in H_0^1(\Omega), R_h\psi \in V_h$ such that
\begin{eqnarray*}
   && (\nabla R_h\psi,\nabla v_h)=(\nabla \psi,\nabla v_h),\qquad \forall v_h\in V_h.
\end{eqnarray*}
The following approximation property is known \cite[Chapter 1]{Thomee-1984}:
\begin{eqnarray*}
   \|R_h\psi-\psi\|_{L^2(\Omega)} + h\|\nabla (R_h\psi-\psi)\|_{L^2(\Omega)}\leq ch^q\|\psi\|_{H^q(\Omega)},\qquad\forall\psi\in H_0^1(\Omega)\cap H^q(\Omega),~q=1,2.
\end{eqnarray*}
Taking the test function in \eqref{lem-stab-ipthm-3-eq1} as the projection of $v\in H_0^1(\Omega)$, i.e.,
$v_{h_m}:=R_{h_m}v$, then
it follows from the above error bound for $q=1$ that $v_{h_m}\to v$ in $L^2(\Omega)$ as $m\to \infty$. Combining this with \eqref{lem-stab-ipthm-3-eq2} and the boundedness of $\|\xi_{h_m}^k\|_{H^1(\Omega)}$ with respect to $h_m$, we arrive at
\begin{eqnarray*}
   && (\xi_{h_m}^k, v_{h_m}) = (  \xi_{h_m}^k, v ) + (  \xi_{h_m}^k, v_{h_m}-v ) \to ( \xi^k, v)\quad\text{as}\quad m\to \infty,
\end{eqnarray*}
for all $k\geq 1$.
Notice that $(\nabla v_{h_m},\nabla \xi_{h_m}^k)=(\nabla v,\nabla \xi_{h_m}^k)$, we have
\begin{eqnarray*}
   && (\nabla \xi_{h_m}^k,\nabla v_{h_m}) = ( \nabla \xi_{h_m}^k,\nabla v ) + ( \nabla \xi_{h_m}^k,\nabla v_{h_m}- \nabla v ) \to (\nabla \xi^k,\nabla v)\quad\text{as}\quad m\to \infty.
\end{eqnarray*}
Therefore, passing $m\to \infty$ in \eqref{lem-stab-ipthm-3-eq1} leads to
\begin{eqnarray}
 \xi^0=0\quad\text{and}\quad(\xi^{k+1},v) + \eta(\nabla \xi^{k+1},\nabla v) = \sum\limits_{j=0}^{k-1}(b_j-b_{j+1}) (\xi^{k-j},v),\quad \forall v\in H_0^1(\Omega). \label{lem-stab-ipthm-3-eq3}
\end{eqnarray}
Obviously this semi-discrete problem admits the unique trivial solution $\xi^k=0, k\geq 0$.
Consequently for all $k\geq 0$ we have $\|\xi_{h_m}^k\|_{L^2(\Omega)}\to 0$ as $m\to \infty$. This completes the proof.
\hfill$\Box$

We now give the convergence of the finite element approximation \eqref{min-L2TV-disc-0}
to the continuous minimization problem \eqref{min-L2TV-contin}.
\begin{thm}\label{thm-convfm-ipthm}
Let $g\in L^2(\Omega)$,
$\{h_m\}_{m=1}^\infty$ and $\{\tau_m\}_{m=1}^\infty$ be any positive sequences such that $h_m\to 0$ and $\tau_m\to 0$ as $m\to \infty$,
$\{f_{m}\}_{m=1}^\infty\subset X_{h_m}^1$ be the minimizers to the discrete optimization problems \eqref{min-L2TV-disc-0} with $h=h_m$ and $\tau=\tau_m$.
Then an element $f^*\in BV(\Omega)\cap L^2(\Omega)$ exists such that
\begin{eqnarray}
   && \lim\limits_{m\to \infty}\|f_{m}-f^*\|_{L^2(\Omega)}=0 \quad\text{and}\quad \lim\limits_{m\to \infty}TV(f_{m})=TV(f^*). \label{thm-convfm-resipthm}
\end{eqnarray}
Furthermore, $f^*$ is the unique solution to the continuous problem \eqref{min-L2TV-contin}.
\end{thm}
\noindent \emph{Proof }
By virtue of Lemma \ref{lem-stab-ipthm-1},
for $f\in BV(\Omega)\cap L^2(\Omega)$ there exists $\widehat{f}_{h_m}\in X_{h_m}^1$ such that
\begin{align}
\lim\limits_{m\to \infty}\|\widehat{f}_{h_m}-f\|_{L^2(\Omega)}=0 \quad \text{and}\quad \lim\limits_{m\to \infty}TV(\widehat{f}_{h_m})=TV(f). \label{thm-convfm-ipthm-eq1}
\end{align}
The minimizers $f_m$ of \eqref{min-L2TV-disc-0} satisfy
\begin{align}
J_{h_m,\tau_m}(f_m)&=\frac{1}{2}\|u_{h_m}^{K_{\tau_m}}(f_m)-g\|_{L^2(\Omega)}^2 + \frac{\beta}{2}\|f_m\|_{L^2(\Omega)}^2 + \gamma TV(f_m) \nonumber\\
&\leq \frac{1}{2}\|u_{h_m}^{K_{\tau_m}}(\widehat{f}_{h_m})-g\|_{L^2(\Omega)}^2 + \frac{\beta}{2}\|\widehat{f}_{h_m}\|_{L^2(\Omega)}^2 + \gamma TV(\widehat{f}_{h_m}). \label{thm-convfm-ipthm-eq2}
\end{align}
Then it follows from \eqref{thm-convfm-ipthm-eq1}, \eqref{thm-convfm-ipthm-eq2}, and Theorem \ref{thm-uncstable-fuly} that $\{f_{m}\}_{m=1}^\infty$
is a bounded set in $BV(\Omega)\cap L^2(\Omega)$. By Proposition \ref{prop-BVfun} $(ii)$, there exists a (not relabeled) subsequence of $\{f_{m}\}_{m=1}^\infty$ and an element $f^*\in BV(\Omega)\cap L^2(\Omega)$ such that
\begin{eqnarray*}
  f_m\rightharpoonup f^*\ \text{in}\ L^2(\Omega), \ \ \
  f_m\to f^*\ \text{in}\ L^1(\Omega),\quad\text{as } m\to \infty,
\end{eqnarray*}
and
\begin{eqnarray}
  TV(f^*)\leq\liminf\limits_{m\to\infty}TV(f_m).\label{thm-convfm-ipthm-eq3}
\end{eqnarray}
Furthermore, using Lemma \ref{lem-stab-ipthm-3} and the lower semi-continuity of the $L^2$-norm, we have
\begin{align*}
J(f^*)&=\frac{1}{2}\|u(f^*)-g\|_{L^2(\Omega)}^2 + \frac{\beta}{2}\|f^*\|_{L^2(\Omega)}^2 + \gamma TV(f^*) \\
&\leq \lim\limits_{m\to \infty}\frac{1}{2}\|u_{h_m}^{K_{\tau_m}}(f_m)-g\|_{L^2(\Omega)}^2 + \liminf\limits_{m\to \infty}\big(\frac{\beta}{2}\|f_m\|_{L^2(\Omega)}^2 + \gamma TV(f_m)\big)\\
&\leq \liminf\limits_{m\to \infty} \big(\frac{1}{2}\|u_{h_m}^{K_{\tau_m}}(f_m)-g\|_{L^2(\Omega)}^2 + \frac{\beta}{2}\|f_m\|_{L^2(\Omega)}^2 + \gamma TV(f_m)\big)\\
& \leq \limsup\limits_{m\to \infty} \big(\frac{1}{2}\|u_{h_m}^{K_{\tau_m}}(f_m)-g\|_{L^2(\Omega)}^2 + \frac{\beta}{2}\|f_m\|_{L^2(\Omega)}^2 + \gamma TV(f_m)\big)\\
&\leq \limsup\limits_{m\to \infty} \big(\frac{1}{2}\|u_{h_m}^{K_{\tau_m}}(\widehat{f}_{h_m})-g\|_{L^2(\Omega)}^2 + \frac{\beta}{2}\|\widehat{f}_{h_m}\|_{L^2(\Omega)}^2 + \gamma TV(\widehat{f}_{h_m})\big)\\
&  =\frac{1}{2}\|u(f)-g\|_{L^2(\Omega)}^2 + \frac{\beta}{2}\|f\|_{L^2(\Omega)}^2 + \gamma TV(f) =J(f).
\end{align*}
This shows that $f^*$ is a solution to the problem \eqref{min-L2TV-contin}.
Now, replacing $f$ by $f^*$ in the above analysis, we can obtain
\begin{align*}
&\frac{1}{2}\|u(f^*)-g\|_{L^2(\Omega)}^2 + \frac{\beta}{2}\|f^*\|_{L^2(\Omega)}^2 + \limsup\limits_{m\to \infty}\gamma TV(f_m)\\
\leq&  \limsup\limits_{m\to \infty} \big(\frac{1}{2}\|u_{h_m}^{K_{\tau_m}}(f_m)-g\|_{L^2(\Omega)}^2 + \frac{\beta}{2}\|f_m\|_{L^2(\Omega)}^2 + \gamma TV(f_m)\big)\\
=&  \frac{1}{2}\|u(f^*)-g\|_{L^2(\Omega)}^2 + \frac{\beta}{2}\|f^*\|_{L^2(\Omega)}^2 + \gamma TV(f^*).
\end{align*}
This together with \eqref{thm-convfm-ipthm-eq3} leads to
\begin{eqnarray*}
  && TV(f^*)\leq\liminf\limits_{m\to \infty}TV(f_m)\leq \limsup\limits_{m\to \infty}TV(f_m)= TV(f^*).
\end{eqnarray*}
Thus $\lim\limits_{m\to \infty}TV(f_m)=TV(f^*)$. Similarly it can be shown $\lim\limits_{m\to \infty}\|f_m\|_{L^2(\Omega)}=\|f^*\|_{L^2(\Omega)}$, which together with the weak convergence implies the strong convergence in $L^2(\Omega)$ of $f_m$ to $f^*$.
Finally, it is not difficult to prove by contradiction that the convergence \eqref{thm-convfm-resipthm} holds for the whole sequence
and the solution of \eqref{min-L2TV-contin} is unique.
\hfill$\Box$

\begin{thm}\label{thm-errfm-ipthm}
Let $g\in L^2(\Omega)$.
If $f_{h,\tau}\in X_h^1(\Omega)$ and $f^*\in BV(\Omega)\cap L^2(\Omega)$ are
the minimizers of
\eqref{min-L2TV-disc-0} and \eqref{min-L2TV-contin}, respectively.
Then
\begin{eqnarray*}
   && \|f_{h,\tau}-f^*\|_{L^2(\Omega)}^2\leq c \left( \frac{\beta+\gamma+c}{\beta}h^{\frac{1}{3}} + \frac{\tau}{\beta} + \frac{h^2|\ln h|^2}{\beta} \right).
\end{eqnarray*}
\end{thm}

\noindent \emph{Proof }
Let $f^\circ\in BV(\Omega)\cap L^\infty(\Omega)$
be the solution of $\min J_{h,\tau}(f)$ over $BV(\Omega)\cap L^2(\Omega)$.
Then it follows from Lemma \ref{lem-stab-ipthm-2}(i) that there exists
$\widehat{f}_h\in X_h^1$ such that
\begin{align}
\|\widehat{f}_{h}-f^\circ\|_{L^2(\Omega)}\leq ch^{\frac{1}{3}} \quad \text{and}\quad TV(\widehat{f}_h)-TV(f^\circ)\leq ch^{\frac{1}{3}}. \label{thm-errfm-ipthm-eq1}
\end{align}
By using Cauchy-Schwarz inequality, the a prior estimates of $u_h^{K_\tau}(\cdot)$, the boundedness of $\{\widehat{f}_h\}_{h>0}$, and \eqref{thm-errfm-ipthm-eq1}, we have
\begin{align*}
  & J_{h,\tau}(\widehat{f}_{h}) - J_{h,\tau}(f^\circ) \\
 = & \frac{1}{2}\|u_h^{K_\tau}(\widehat{f}_{h})-g\|_{L^2(\Omega)}^2 +\frac{\beta}{2}\|\widehat{f}_h\|_{L^2(\Omega)}^2 +\gamma TV(\widehat{f}_h) -\frac{1}{2}\|u_h^{K_\tau}(f^\circ)-g\|_{L^2(\Omega)}^2 -\frac{\beta}{2}\|f^\circ\|_{L^2(\Omega)}^2 -\gamma TV(f^\circ)\\
 = & \frac{\beta}{2}(\widehat{f}_h-f^\circ,\widehat{f}_h+f^\circ) +\gamma(TV(\widehat{f}_h)-TV(f^\circ)) +\frac{1}{2}(u_h^{K_\tau}(\widehat{f}_{h})+u_h^{K_\tau}(f^\circ)-2g,u_h^{K_\tau}(\widehat{f}_{h})-u_h^{K_\tau}(f^\circ)) \\
 \leq & c\beta \|\widehat{f}_h-f^\circ\|_{L^2(\Omega)} +c\gamma h^{\frac{1}{3}} +c\|u_h^{K_\tau}(\widehat{f}_{h})-u_h^{K_\tau}(f^\circ)\|_{L^2(\Omega)} \\
 \leq & c(\beta+1) \|\widehat{f}_h-f^\circ\|_{L^2(\Omega)} +c\gamma h^{\frac{1}{3}} \\
 \leq & c(\beta+\gamma+c)h^{\frac{1}{3}}.
\end{align*}
Furthermore, it follows from the $\beta$-strong convexity of $J_{h,\tau}(\cdot)$
and the optimality of $f^\circ$ and $f_{h,\tau}$:
\begin{eqnarray}
   \frac{\beta}{2}\|f_{h,\tau}-f^\circ\|_{L^2(\Omega)}^2\leq J_{h,\tau}(f_{h,\tau})-J_{h,\tau}(f^\circ)\leq  J_{h,\tau}(\widehat{f}_h)-J_{h,\tau}(f^\circ)\leq c(\beta+\gamma+c)h^{\frac{1}{3}}.\label{thm-errfm-ipthm-eq2}
\end{eqnarray}
On the other hand, by Theorem \ref{thm-fuly-erest}, we have
\begin{align*}
J(f^\circ) &= \frac{1}{2}\|u(f^\circ)-g\|_{L^2(\Omega)}^2 +\frac{\beta}{2}\|f^\circ\|_{L^2(\Omega)}^2 +\gamma TV(f^\circ) \\
& =\frac{1}{2}\|u(f^\circ)-u_h^{K_\tau}(f^\circ)\|_{L^2(\Omega)}^2 + (u(f^\circ)-u_h^{K_\tau}(f^\circ),u_h^{K_\tau}(f^\circ)-g) \\
& \quad +\frac{1}{2}\|u_h^{K_\tau}(f^\circ)-g\|_{L^2(\Omega)}^2 +\frac{\beta}{2}\|f^\circ\|_{L^2(\Omega)}^2 +\gamma TV(f^\circ) \\
& \leq c\|u(f^\circ)-u_h^{K_\tau}(f^\circ)\|_{L^2(\Omega)} + J_{h,\tau}(f^\circ) \\
& \leq c(\tau + h^2|\ln h|^2) + J_{h,\tau}(f^\circ).
\end{align*}
Again, the strong convexity of $J(\cdot)$ leads to
\begin{align*}
\frac{\beta}{2}\|f^\circ-f^*\|_{L^2(\Omega)}^2 & \leq J(f^\circ)-J(f^*) \\
& \leq J_{h,\tau}(f^\circ) -J(f^*) + c(\tau + h^2|\ln h|^2) \\
& \leq J_{h,\tau}(f^*) -J(f^*) + c(\tau + h^2|\ln h|^2) \\
& =\frac{1}{2}\|u_h^{K_\tau}(f^*)-g\|_{L^2(\Omega)}^2 -\frac{1}{2}\|u(f^*)-g\|_{L^2(\Omega)}^2 +c(\tau + h^2|\ln h|^2) \\
& = \frac{1}{2}(u_h^{K_\tau}(f^*)-u(f^*), u_h^{K_\tau}(f^*)+u(f^*)-2g) +c(\tau + h^2|\ln h|^2)\\
& \leq c\|u_h^{K_\tau}(f^*)-u(f^*)\|_{L^2(\Omega)} +c(\tau + h^2|\ln h|^2)\\
& \leq c(\tau + h^2|\ln h|^2).
\end{align*}
Finally we combine the above estimate with \eqref{thm-errfm-ipthm-eq2}
to conclude.
\hfill$\Box$

\begin{rem}
If we suppose that $f^\circ\in BV(\Omega)\cap \mathrm{Lip}(s,L^2(\Omega))$ with $0<s\leq 1$, then a similar analysis will give
\begin{eqnarray*}
   && \|f_{h,\tau}-f^*\|_{L^2(\Omega)}^2\leq c \left( \frac{\beta+\gamma+c}{\beta}h^{\frac{s}{s+1}} + \frac{\tau}{\beta}+ \frac{h^2|\ln h|^2}{\beta} \right).
\end{eqnarray*}
\end{rem}

In the remaining part of this section, we analyze the stability of the solution to
\eqref{min-L2TV-disc-0} with respect to the data $g$. To this end, we consider
\begin{align}
 \min\limits_{f\in \mathcal{I}}  \frac{\beta}{2}\|f\|_{L^2(\Omega)}^2 + \gamma TV(f), \quad\text{with}\quad \mathcal{I}:=\{f\in BV(\Omega)\cap L^2(\Omega)~|~ u(f)=g\},   \label{min-idenp}
\end{align}
for some $\beta,\gamma>0$. The solution of \eqref{min-idenp} can be regarded
as a ``minimal" solution to the unregularized least squares problem
\begin{align*}
 \min\limits_{f\in BV(\Omega)\cap L^2(\Omega)} \frac{1}{2}\|u(f)-g\|_{L^2(\Omega)}^2.
\end{align*}
We know from standard arguments of convex analysis that the problem \eqref{min-idenp} admits a unique solution.
Next we give the convergence results to a sought source functional when the regularization parameter approaches to zero with a suitable coupling of noise level and mesh size.

\begin{thm}\label{thm-stab-data}
Assume that $f^\dag\in BV(\Omega)\cap L^\infty(\Omega)$ is the unique solution of \eqref{min-idenp}.
Let $\{\delta_m\}_{m=1}^\infty$, $\{h_m\}_{m=1}^\infty$, $\{\tau_m\}_{m=1}^\infty$, $\{\beta_m\}_{m=1}^\infty$, and $\{\gamma_m\}_{m=1}^\infty$ be positive zero sequences as $m\to \infty$,
and $\{g^{\delta_m}\}_{m=1}^\infty\subset L^2(\Omega)$ be a sequence of noisy data satisfying $\|g^{\delta_m}-g\|_{L^2(\Omega)}\leq \delta_m$.
Let $\{f^{\delta_m}\}_{m=1}^\infty$ be the minimizers to the discrete optimization problems \eqref{min-L2TV-disc-0}
associated to the noisy data $g^{\delta_m}$
with $h=h_m$, $\tau=\tau_m$, $\beta=\beta_m$, and $\gamma=\gamma_m$.
Then
\begin{eqnarray}
   &&  \lim\limits_{m\to\infty}\|f^{\delta_m} -f^\dag\|_{L^2(\Omega)}=0 \quad\text{and}\quad \lim\limits_{m\to\infty}TV(f^{\delta_m})=TV(f^\dag)
\label{lem-stab-ipthm-0}
\end{eqnarray}
when
\begin{eqnarray}
 \beta_m\to 0,\quad \frac{\beta_m}{\gamma_m}\to\frac{\beta}{\gamma},\quad \frac{\delta_m}{\sqrt{\beta_m}}\to 0,
 \quad \frac{h_m|\ln h_m|}{\beta_m}\to 0, \quad\text{and}\quad \frac{\tau_m}{\sqrt{\beta_m}}\to 0 \quad\text{as}\quad m\to\infty. \label{thm-cond-data}
\end{eqnarray}
\end{thm}

\noindent \emph{Proof }
It follows from the optimality of $f^{\delta_m}$ that
\begin{align*}
J_{h_m,\tau_m}^{\delta_m}(f^{\delta_m})&=\frac{1}{2}\|u_{h_m}^{K_{\tau_m}}(f^{\delta_m})-g^{\delta_m}\|_{L^2(\Omega)}^2 + \frac{\beta_m}{2}\|f^{\delta_m}\|_{L^2(\Omega)}^2 + \gamma_m TV(f^{\delta_m}) \nonumber\\
&\leq \frac{1}{2}\|u_{h_m}^{K_{\tau_m}}(\widehat{f}_{h_m})-g^{\delta_m}\|_{L^2(\Omega)}^2 + \frac{\beta_m}{2}\|\widehat{f}_{h_m}\|_{L^2(\Omega)}^2 + \gamma_m TV(\widehat{f}_{h_m}),
\end{align*}
According to Lemma \ref{lem-stab-ipthm-2}(ii), for $f^\dag\in\mathcal{I}$ there exist $\widehat{f}_{h_m}\in X_{h_m}^1$ such that
\begin{align}
\|\widehat{f}_{h_m}-f^\dag\|_{L^p(\Omega)}\leq c(h_m|\ln h_m|)^{\frac{1}{p}}
\quad \text{and}\quad \lim\limits_{h_m\to 0}TV(\widehat{f}_{h_m})=TV(f^\dag).
\label{unapBV-inf-1b2}
\end{align}
Then
\begin{align*}
    \frac{1}{2}\|u_{h_m}^{K_{\tau_m}}(\widehat{f}_{h_m})-g\|_{L^2(\Omega)}^2
    &\leq \|u_{h_m}^{K_{\tau_m}}(\widehat{f}_{h_m})-u(f^\dag)\|_{L^2(\Omega)}^2 + \|u(f^\dag)-g\|_{L^2(\Omega)}^2 \\
    &\leq \|u_{h_m}^{K_{\tau_m}}(\widehat{f}_{h_m})-u(f^\dag)\|_{L^2(\Omega)}^2 +\delta_m^2.
\end{align*}
Furthermore, applying Theorem \ref{thm-uncstable-fuly}, Lemma \ref{lem-stab-ipthm-2}, and Theorem \ref{thm-fuly-erest},
the first term on the right side can be bounded by
\begin{align*}
     \|u_{h_m}^{K_{\tau_m}}(\widehat{f}_{h_m})-u(f^\dag)\|_{L^2(\Omega)} & \leq \|u_{h_m}^{K_{\tau_m}}(\widehat{f}_{h_m})-u_{h_m}^{K_{\tau_m}}(f^\dag)\|_{L^2(\Omega)} + \|u_{h_m}^{K_{\tau_m}}(f^\dag)-u(f^\dag)\|_{L^2(\Omega)} \\
    &\leq c\|\widehat{f}_{h_m}-f^\dag\|_{L^2(\Omega)} +  \|u_{h_m}^{K_{\tau_m}}(f^\dag)-u(f^\dag)\|_{L^2(\Omega)} \\
    & \leq c(h_m|\ln h_m|)^{\frac{1}{2}} + c(\tau_m +h_m^2|\ln h_m|^2) \\
    &\leq c((h_m|\ln h_m|)^{\frac{1}{2}}+\tau_m).
\end{align*}
Combining the above results gives
\begin{align}
    & \frac{1}{2}\|u_{h_m}^{K_{\tau_m}}(f^{\delta_m})-g^{\delta_m}\|_{L^2(\Omega)}^2 + \frac{\beta_m}{2}\|f^{\delta_m}\|_{L^2(\Omega)}^2 + \gamma_m TV(f^{\delta_m})  \nonumber\\
    \leq & c( h_m|\ln h_m|+\tau_m^2 + \delta_m^2 ) + \frac{\beta_m}{2}\|\widehat{f}_{h_m}\|_{L^2(\Omega)}^2 + \gamma_m TV(\widehat{f}_{h_m}). \label{thm-stab-data-eq1}
\end{align}
Therefore, by using \eqref{thm-cond-data} and  Lemma \ref{lem-stab-ipthm-2}, we have
\begin{eqnarray}
   &&  \lim\limits_{m\to\infty} \|u_{h_m}^{K_{\tau_m}}(f^{\delta_m})-g^{\delta_m}\|_{L^2(\Omega)}=0, \label{thm-stab-data-eq2}\\
   && \limsup\limits_{m\to\infty} \frac{1}{2}\|f^{\delta_m}\|_{L^2(\Omega)}^2 + \frac{\beta_m}{\gamma_m} TV(f^{\delta_m})\leq
   \frac{1}{2}\|f^\dag\|_{L^2(\Omega)}^2 + \frac{\beta}{\gamma} TV(f^\dag). \label{thm-stab-data-eq3}
\end{eqnarray}
It is routine to get from \eqref{thm-stab-data-eq1}
the boundedness of $\{f^{\delta_m}\}_{m=1}^\infty\subset BV(\Omega)\cap L^2(\Omega)$.
In virtue of Proposition \ref{prop-BVfun}(ii) and Lemma \ref{lem-stab-ipthm-3}, a (not relabeled) subsequence of
$\{f^{\delta_m}\}_{m=1}^\infty$  and an element $f^\diamond\in BV(\Omega)\cap L^2(\Omega)$ exist such that
\begin{eqnarray*}
  & f^{\delta_m}\rightharpoonup f^\diamond\ \text{in}\ L^2(\Omega),\quad f^{\delta_m}\to f^\diamond\ \text{in}\ L^1(\Omega),  \quad\text{as } m\to \infty, & \\
  & TV(f^\diamond)\leq\liminf\limits_{m\to\infty}TV(f^{\delta_m}), &
\end{eqnarray*}
and
\begin{eqnarray}
  \lim\limits_{m\to\infty}\|u_{h_m}^{K_{\tau_m}}(f^{\delta_m})-g\|_{L^2(\Omega)}=\|u(f^\diamond)-g \|_{L^2(\Omega)}. \label{thm-stab-data-eq4}
\end{eqnarray}
Then it follows from \eqref{thm-stab-data-eq2} and \eqref{thm-stab-data-eq4} that
\begin{eqnarray*}
\|u(f^\diamond)-g\|_{L^2(\Omega)}\leq \lim\limits_{m\to\infty}\big(  \|u(f^\diamond)-u_{h_m}^{K_{\tau_m}}(f^{\delta_m})\|_{L^2(\Omega)}
 + \|u_{h_m}^{K_{\tau_m}}(f^{\delta_m})-g^{\delta_m}\|_{L^2(\Omega)} + \|g^{\delta_m}-g\|_{L^2(\Omega)} \big)=0.
\end{eqnarray*}
This implies $f^\diamond\in\mathcal{I}$. Furthermore, it follows from \eqref{thm-stab-data-eq3}
\begin{align*}
\frac{1}{2}\|f^\diamond\|_{L^2(\Omega)}^2 + \frac{\beta}{\gamma} TV(f^\diamond) &\leq \liminf\limits_{m\to\infty}\big( \frac{1}{2}\|f^{\delta_m}\|_{L^2(\Omega)}^2 + \frac{\beta_m}{\gamma_m} TV(f^{\delta_m}) \big) \\
&\leq \limsup\limits_{m\to\infty}\big( \frac{1}{2}\|f^{\delta_m}\|_{L^2(\Omega)}^2 + \frac{\beta_m}{\gamma_m} TV(f^{\delta_m}) \big)\\
&\leq \frac{1}{2}\|f^\dag\|_{L^2(\Omega)}^2 + \frac{\beta}{\gamma} TV(f^\dag).
\end{align*}
Hence $f^\diamond=f^\dag$ by the uniqueness of the solution of \eqref{min-idenp}. Finally,
similar to the proof of Theorem \ref{thm-convfm-ipthm}, we obtain \eqref{lem-stab-ipthm-0}.
This ends the proof.
\hfill$\Box$

\begin{rem}
If $f^\dag\in\mathrm{Lip}(s,L^2(\Omega)), 0<s\leq 1$, then the convergence condition
of \eqref{lem-stab-ipthm-0} on the parameters would become
\begin{eqnarray*}
 \beta_m\to 0,\quad \frac{\beta_m}{\gamma_m}\to\frac{\beta}{\gamma},\quad \frac{\delta_m}{\sqrt{\beta_m}}\to 0,\quad \frac{h_m|\ln h_m|}{\sqrt[2s]{\beta_m}}\to 0\quad\text{and}\quad \frac{\tau_m}{\sqrt{\beta_m}}\to 0 \quad\text{as}\quad m\to\infty.
\end{eqnarray*}
\end{rem}

\section{An accelerated linearized primal-dual algorithm}\label{sec:4}

In this section, we propose an algorithm to find minimizers of the problem \eqref{min-L2TV-disc-0} and
carry out error analysis of the proposed algorithm.

\subsection{ Saddle point problem and optimality}

The non-smooth TV items in the minimization problems present more challenges
in designing efficient iterative schemes.
Our algorithm is based on the TV-dual representation approach
\cite{Bartels-2012,Bartels-2015,Chambolle-2011,Tian-Yuan-2019}, which consists in reformulating
the minimization (\ref{min-L2TV-contin}) into a saddle point problem as follows:
\begin{eqnarray}\label{sadd-point-prob}
\inf\limits_{f\in BV(\Omega)\cap L^2(\Omega)}\sup\limits_{\rho\in \mathcal{B}} \Psi^\delta(f,\rho),
\end{eqnarray}
where
\begin{eqnarray*}
\Psi^\delta(f,\rho):=\frac{1}{2}\|u(f)-g^\delta\|_{L^2(\Omega)}^2 + \frac{\beta}{2}\|f\|_{L^2(\Omega)}^2
+ \gamma\int_\Omega f\nabla\cdot\rho ~\mathrm{d}x.
\end{eqnarray*}
The approximation to the vector function $\rho$ makes use of the piecewise constant finite element space $X_h^0$,
defined by
\begin{eqnarray*}
X_h^0 = \left\{q_h\in L^1(\Omega)~ | ~ q_h|_K=\text{constant},~\forall K\in\mathcal{T}_h \right\}.
\end{eqnarray*}
Noticing that \cite{Bartels-2012}
\begin{eqnarray}
 TV(f)=\int_{\Omega}|\nabla f|\mathrm{d}x=\sup\limits_{\rho\in \mathcal{B}_1}
 \int_{\Omega}\nabla f\cdot\rho~\mathrm{d} x,\qquad \forall f\in X_h^1,
 \label{secdef-TV}
\end{eqnarray}
where $\mathcal{B}_1:=\{\rho\in(X_h^0)^d~|~|\rho|_\infty\leq 1 \}$.
We propose to approximate the saddle point problem \eqref{sadd-point-prob} by a mixed finite element method as follows:
\begin{eqnarray}\label{dis-sadd-point-prob}
 \inf\limits_{f\in X_h^1}\sup\limits_{\rho\in (X_h^0)^d} \Psi^\delta_{h,\tau}(f,\rho),
\end{eqnarray}
where
\begin{eqnarray*}
\Psi_{h,\tau}^\delta(f,\rho):=\frac{1}{2}\|u_h^{K_\tau}(f)-g^\delta\|_{L^2(\Omega)}^2 + \frac{\beta}{2}\|f\|_{L^2(\Omega)}^2 + \gamma\int_{\Omega}\nabla f\cdot \rho\mathrm{d}x- \delta_{\mathcal{B}_1}(\rho),
\end{eqnarray*}
$\delta_{\mathcal{B}_1}(\cdot)$ denotes the indicator function of the set $\mathcal{B}_1$.

Denote $(X_h^1)^*$ and $((X_h^0)^d)^*$ the dual spaces of $X_h^1$ and $(X_h^0)^d$, respectively.
We first list a number of known results in the following remark.

\begin{rem}\label{rem-der-dual}\

$1^\circ$~ It is readily seen that $u_h^{K_\tau}(f)$ is linear with respect to $f$, and the Fr\'{e}chet derivative
$u_h^{K_\tau}(f)'z=u_h^{K_\tau}(z)$ for any $z\in X_h^1$.

$2^\circ$~Define the norm of the weak gradient operator $\nabla$ by
\begin{eqnarray*}
   && \|\nabla\|:=\sup_{0\neq f\in X_h^1} \frac{\|\nabla f\|_{L^2(\Omega)}}{\|f\|_{L^2(\Omega)}}.
\end{eqnarray*}
Then it follows from the well known inverse inequality in $X_h^1$ that $\|\nabla\|\leq ch^{-1}$.

$3^\circ$~ The inclusions $X_h^1\subset (X_h^1)^*$ and $(X_h^0)^d\subset ((X_h^0)^d)^*$ hold via the identities
\begin{eqnarray*}
   &&  (z,f)_{((X_h^1)^*,X_h^1)} := (z,f),\qquad \forall f\in X_h^1,\quad\forall z\in X_h^1, \\
   && (q,\rho)_{ (((X_h^0)^d)^*,(X_h^0)^d) } := (q,\rho),\qquad \forall \rho\in (X_h^0)^d,\quad\forall q\in (X_h^0)^d.
\end{eqnarray*}

$4^\circ$~Any $\rho\in (X_h^0)^d$ can be considered as an element of $(X_h^1)^*$, by
\begin{eqnarray*}
   && (\rho,f)_{((X_h^1)^*,X_h^1)} := (\nabla f,\rho),\qquad\forall f\in X_h^1,\quad\forall \rho\in(X_h^0)^d.
\end{eqnarray*}

$5^\circ$~ For each $f\in X_h^1$ the relation
\begin{eqnarray}
   && \partial TV(f)=\left\{ \rho\in\mathcal{A}\subset(X_h^0)^d\subset (X_h^1)^*~|~(\nabla f,\rho)=\int_{\Omega}|\nabla f|\mathrm{d}x\right\} \label{subgr-TV}
\end{eqnarray}
holds, where $\partial TV(f)$ is the subgradient of $f$; see \cite[Lemma 4.2]{Hinze-Quyen-2019} for a detailed proof.
\end{rem}

The next lemma gives the first-order optimality condition of the problem \eqref{dis-sadd-point-prob}.

\begin{lem}\label{lem-optimality}$($Optimality$)$
The function $f\in X_h^1$ is a solution of \eqref{min-L2TV-disc-0} if and only if there exists $\rho\in \partial TV(f)$ such that
\begin{eqnarray}
   &&  (u_h^{K_\tau}(f)-g^\delta),u_h^{K_\tau}(z)) + \beta(f,z) + \gamma(\rho,\nabla z)=0, \qquad \forall z\in X_h^1, \label{lem-opti-cond1}\\
   && (\nabla f,q-\rho)\leq 0, \qquad  \forall q\in\mathcal{A}.  \label{lem-opti-cond2}
\end{eqnarray}
\end{lem}
\noindent \emph{Proof }
The optimality condition comes from the standard Kuhn-Tucker conditions, i.e.,
\begin{eqnarray}
   && 0=\partial_f \Psi_{h,\tau}^\delta(f,\rho),\qquad 0\in \partial_\rho \Psi_{h,\tau}^\delta(f,\rho). \label{lem-optimality-eq1}
\end{eqnarray}
According to Remark \ref{rem-der-dual} $1^\circ$ and $4^\circ$,
we obtain \eqref{lem-opti-cond1} from the first equation of \eqref{lem-optimality-eq1}.
Furthermore, there exists $\rho\in\mathcal{A}$ such that $\nabla f\in\partial\delta_{\mathcal{A}}(\rho)$,
where the subgradient of $\delta_{\mathcal{A}}(\cdot)$ is given by \cite[Example 3.19]{Peypouquet-2015}:
\begin{eqnarray*}
   &&  \partial\delta_{\mathcal{A}}(\rho)=\left\{ \varrho\in((X_h^0)^d)^* ~|~ (\varrho,q-\rho) \leq 0 \quad\text{for~all }q\in\mathcal{A} \right\}.
\end{eqnarray*}
Thus we have \eqref{lem-opti-cond2}. Finally, \eqref{lem-opti-cond2} implies
$(\nabla f,\rho)\geq (\nabla f,q)$ for all $q\in\mathcal{A}$. This, together with \eqref{secdef-TV}, gives
$(\nabla f,\rho)=TV(f)$. Hence it follows from Remark \ref{rem-der-dual} $5^\circ$ that
$\rho\in\partial TV(f)$.
\hfill$\Box$

\begin{rem}\label{rem-optimality}
We can rewrite the system \eqref{lem-opti-cond1}-\eqref{lem-opti-cond2} as the following variational inequality in a compact form:
find $\nu\in X_h^1\times\mathcal{A}$, such that
\begin{eqnarray}
   && (F(\nu),\kappa-\nu)\geq 0, \qquad \forall \kappa\in X_h^1\times\mathcal{A}, \label{rem-optimality-comp}
\end{eqnarray}
where
\begin{eqnarray*}
  \nu:=\left(\begin{array}{c}
               f \\
               \rho
             \end{array}\right)\in X_h^1\times\partial TV(f)\subset X_h^1\times\mathcal{A},\qquad \kappa:=\left(\begin{array}{c}
               z \\
               q
             \end{array}\right)\in X_h^1\times\mathcal{A},
\end{eqnarray*}
and
\begin{eqnarray*}
   && F(\nu):=\left(\begin{array}{c}
                      (u_h^{K_\tau})^*(u_h^{K_\tau}(f)-g^\delta) + \beta f-\gamma\mathrm{div}\rho \\
                      -\beta\nabla f
                    \end{array}\right),
\end{eqnarray*}
in which $(u_h^{K_\tau})^*$ and $-\mathrm{div}$ denote the adjoint operators of $u_h^{K_\tau}$ and $\nabla$, respectively.

\end{rem}

\subsection{Iterative algorithm}
Now we propose a specific iterative method to solve the saddle point problem \eqref{dis-sadd-point-prob}.
It has been known that for an ill-posed problem, if the solution does not have a certain degree of regularity, the convergence of general solution methods will be very slow;
see, e.g., \cite[Chapter 3]{Engl-Hanke-Beyvayer-1996} for more details.
Based on an idea in \cite{Tian-Yuan-2019},
we adopt here an accelerated linearized primal-dual algorithm to solve
\eqref{dis-sadd-point-prob}.

\begin{alg} (Accelerated linearized primal-dual algorithm)
\label{alg-acc-pd}

\noindent - Input: Let parameters $\varsigma_0$, $\upsilon_0$, $\theta_0>0$ such that
\begin{eqnarray}
   && \frac{1-3c^2\varsigma_0}{\varsigma_0} >\beta^2\frac{\varsigma_0}{\upsilon_0}\|\nabla\|^2, \label{cond1-AG-conv}
\end{eqnarray}
where $c$ is the norm of the operator $u_h^{K_\tau}(\cdot)$. Choose an initial guess
$(f^0,\rho^0)\in X_h^1\times (X_h^0)^d$, $f^{-1}:=f^0$.

\noindent - For $n=0,1,2\cdots,$ do

Update the new iteration $(f^{n+1},\rho^{n+1})$ via solving
 \begin{subequations}\label{Ag-iter-pd}
\begin{align}
    \widetilde{f}^n:&=f^n+\theta_n(f^n-f^{n-1}),  \label{iter-pd-a} \\
   \rho^{n+1}:&=\mathop{\arg\max}\limits_{\rho\in (X_h^0)^d}\left\{ \gamma\int_{\Omega}\nabla \widetilde{f}^n\cdot\rho^n\mathrm{d}x -\delta_{\mathcal{A}}(\rho^n) -\frac{\upsilon_n}{2\varsigma_n}\|\rho-\rho^n\|_{L^2(\Omega)}^2  \right\}, \label{iter-pd-b}  \\
    f^{n+1}:&=\mathop{\arg\min}\limits_{f\in X_h^1}\left\{ (u_h^{K_\tau}(f^n)-g^\delta,u_h^{K_\tau}(f))+\frac{\beta}{2}\|f\|_{L^2(\Omega)}^2 +\gamma \int_{\Omega}\nabla f\cdot\rho^{n+1}\mathrm{d}x + \frac{1}{2\varsigma_n}\|f-f^n\|_{L^2(\Omega)}^2\right\}, \label{iter-pd-c}
\end{align}
\end{subequations}

where $\varsigma_n,\upsilon_n$, and $\theta_n$ are updated by
\begin{eqnarray*}
   && \theta_{n+1}:=\frac{1}{\sqrt{1+2\beta\varsigma_n}},\quad \varsigma_{n+1}:=\theta_{n+1}\varsigma_n,\quad \frac{\upsilon_n}{\varsigma_n}=\frac{\upsilon_{n+1}}{\theta_{n+1}\varsigma_{n+1}}.
\end{eqnarray*}

\noindent \ \ end

\noindent - Output: An approximation solution $f^n$ if $n\leq N_{\max}$ for some $N_{\max}\in\mathds{N}$.

\end{alg}

\begin{rem}\

$(1^\circ)$~\eqref{iter-pd-a} accelerates the algorithm to some extent, similar to the idea of Nesterov's strategy \cite{Nesterov-1983}.

$(2^\circ)$~The optimality of $f^{n+1}$ and $\rho^{n+1}$ in the step \eqref{Ag-iter-pd}
yields
\begin{eqnarray}
   && \big( -\frac{\upsilon_n}{\varsigma_n}(\rho^{n+1}-\rho^{n} )+\gamma \nabla \widetilde{f}^n,q-\rho^{n+1} \big)\leq 0,\quad \forall q\in\mathcal{A},  \label{optimality-PD-rho}\\
   && \big( \frac{1}{\varsigma_n}(f^{n+1}-f^n)+\beta f^{n+1},z\big )+(u_h^{K_\tau}(f^n)-g^\delta,u_h^{K_\tau}(z))+\gamma(\rho^{n+1},\nabla z)=0,\quad\forall z\in X_h^1,   \label{optimality-PD-f}
\end{eqnarray}
which can be viewed as a discretization of the $L^2$ gradient flow of the problem \eqref{dis-sadd-point-prob}:
\begin{eqnarray*}
   && -\upsilon \frac{\partial\rho}{\partial t}+\gamma\nabla f\in \partial\delta_{\mathcal{A}}(\rho)
\end{eqnarray*}
and
\begin{eqnarray*}
   && \frac{\partial f}{\partial t} +(u_h^{K_\tau})^*(u_h^{K_\tau}(f)-g^\delta)+\beta f-\gamma\mathrm{div}\rho=0.
\end{eqnarray*}

$(3^\circ)$~The solution $\rho^{n+1}$ to the subproblem \eqref{iter-pd-b} is explicitly given by
\begin{eqnarray}
   && \rho^{n+1}=\frac{\rho^n+\frac{\gamma\varsigma_n}{\upsilon_n}\nabla\widetilde{f}^n }{\max\big\{ 1,|\rho^n+\frac{\gamma\varsigma_n}{\upsilon_n}\nabla\widetilde{f}^n| \big\}},  \label{exfor-PD-rho}
\end{eqnarray}
which can be computed element-wise (refer to \cite{Bartels-2012,Tian-Yuan-2016,Tian-Yuan-2019}).

\end{rem}

Finally, a convergence result of the accelerated linearized primal-dual algorithm
is given without proof
in the following theorem.
A proof can be given by following the lines in
\cite[Theorem 4.1 and Theorem 5.3]{Tian-Yuan-2019}.

\begin{thm}\label{thm-PD-converge}
Let $\{(f^n,\rho^n)\}_{n\in \mathds{N}}$ be the sequence generated by Algorithm \ref{alg-acc-pd}.

$(i)$~Then an element $(f^*,\rho^*)\in X_h^1\times \partial TV(f^*)$ exists such that $\{(f^n,\rho^n)\}_{n\in \mathds{N}}$ converges to $(f^*,\rho^*)$ in $X_h^1\times (X_h^0)^d$. Furthermore, $(f^*,\rho^*)$ satisfies \eqref{rem-optimality-comp} and $f^*$ is a minimizer of \eqref{min-L2TV-disc-0}.

$(ii)$~There holds
\begin{eqnarray*}
   && \|f^{n+1}-f^n\|_{L^2(\Omega)}^2= O\left(\frac{1}{n^2}\right).
\end{eqnarray*}

\end{thm}

\section{Numerical validation} \label{sec:5}

\subsection{Implementation}

Note that $V_h\subset X_h^1$, then the discretized information of the optimization problem \eqref{dis-sadd-point-prob}
can be fully utilized in the calculation of the direct problem \eqref{dir-prob}.
In Algorithm \ref{alg-acc-pd}, the iterations $f^n$ and $\rho^n$ can be computed by \eqref{optimality-PD-f} and \eqref{exfor-PD-rho}, respectively.
However, the term $(u_h^{K_\tau}(f^n)-g^\delta),u_h^{K_\tau}(z))$ in \eqref{optimality-PD-f} implies that we need to solve
\eqref{fuly-dis-dirp} for every direction $z\in X_h^1$.
Direct computation may be acceptable for one-dimensional problems. However,
As it is known, the memory effect of the non-local fractional derivative makes both the computation and memory expensive,
direct calculation can be problematic in the cases of long time integration or high dimensional space.
Here, as an alternative to the direct calculation, we propose the adjoint approach.

Define the functional $\mathcal{H}(f):=\frac{1}{2}\|u(f)-g^\delta\|_{L^2(\Omega)}^2$ for any $f\in L^2(\Omega)$. Its gradient $\mathcal{H}'(f)$ is given by $\mathcal{H}'(f)z=(u(f)-g^\delta,u(z))$ for any $z\in L^2(\Omega)$.
Consider the following adjoint problem
\begin{eqnarray}
   && \left\{\begin{array}{ll}
               D_{T^-}^\alpha w(x,t)-\Delta w(x,t)=0, & x\in\Omega,\quad t\in(0,T], \\
               I_{T^-}^{1-\alpha}w(x,T)=u(f)(x)-g^\delta(x),& x\in\Omega, \\
               w(x,t)=0, &  x\in\partial\Omega,\quad t\in[0,T],
             \end{array}\right. \label{adj-prob}
\end{eqnarray}
where $ D_{T^-}^\alpha$ and $ I_{T^-}^\alpha$ are the right-sided Riemann-Liouville fractional derivative and integral, respectively, defined by
\begin{eqnarray*}
   &&  D_{T^-}^\alpha w(x,t):=-\frac{1}{\Gamma(1-\alpha)} \frac{\mathrm{d}}{\mathrm{d}t} \int_t^T \frac{w(x,s)}{(s-t)^\alpha}\mathrm{d}s,\qquad 0<\alpha<1,\quad 0<t\leq T, \\
   && I_{T^-}^\alpha w(x,t):=\frac{1}{\Gamma(\alpha)}  \int_t^T \frac{w(x,s)}{(s-t)^{1-\alpha}}\mathrm{d}s,\qquad 0<\alpha<1,\quad 0<t\leq T.
\end{eqnarray*}
Then $\mathcal{H}'(f)z$ can be expressed as
\begin{eqnarray}
   && \mathcal{H}'(f)z=\int_\Omega\int_0^T\mu(t)w(x,t)z\mathrm{d} t\mathrm{d} x,  \qquad \forall z\in L^2(\Omega), \label{grad-mH}
\end{eqnarray}
where $w(x,t)$ is solution of \eqref{adj-prob}.
In order to solve \eqref{adj-prob} explicitly, let $\widetilde{w}(x,t)=w(x,T-t)$, then we are led to solve
\begin{eqnarray}
   && \left\{\begin{array}{ll}
               D_{0^+}^\alpha \widetilde{w}(x,t)-\Delta \widetilde{w}(x,t)=0, & x\in\Omega,\quad t\in(0,T], \\
               \lim\limits_{t\to 0^+}I_{0^+}^{1-\alpha}\widetilde{w}(x,t)=u(f)(x)-g^\delta(x),& x\in\Omega, \\
               \widetilde{w}(x,t)=0, &  x\in\partial\Omega,\quad t\in[0,T],
             \end{array}\right. \label{adj-prob-2}
\end{eqnarray}
where $ D_{0^+}^\alpha$ and $ I_{0^+}^\alpha$ are the left-sided fractional derivative and integral:
\begin{eqnarray}
   &&  D_{0^+}^\alpha \widetilde{w}(x,t):=\frac{1}{\Gamma(1-\alpha)} \frac{\mathrm{d}}{\mathrm{d}t} \int_0^t \frac{\widetilde{w}(x,s)}{(t-s)^\alpha}\mathrm{d}s,\qquad 0<\alpha<1,\quad 0<t\leq T, \label{def-RLD-left}\\
   && I_{0^+}^\alpha \widetilde{w}(x,t):=\frac{1}{\Gamma(\alpha)}  \int_0^t \frac{\widetilde{w}(x,s)}{(t-s)^{1-\alpha}}\mathrm{d}s,\qquad 0<\alpha<1,\quad 0<t\leq T.\nonumber
\end{eqnarray}
The solution to \eqref{adj-prob-2} is given by \cite[Section 3]{Yan-Wei-2019}:
\begin{eqnarray}
   && \widetilde{w}(x,t)=\sum\limits_{j=1}^\infty(u(f)(x)-g^\delta(x),\varphi_j)t^{\alpha-1}E_{\alpha,\alpha}(-\lambda_jt^\alpha)\varphi_j \label{formsol-adjpb}
\end{eqnarray}
with $(\lambda_j,\varphi_j)$ being defined by \eqref{varphi-n} and \eqref{lambda-n}.

In practice, the adjoint problem \eqref{adj-prob-2} can be solved in a same way as the direct problem \eqref{dir-prob}.
However more efficient way exists, which is described below.
Using the link between the Riemann-Liouville and Caputo fractional derivatives:
\begin{eqnarray*}
   && D_{0^+}^\alpha \widetilde{w}(x,t)=\partial_t^\alpha \widetilde{w}(x,t) + \frac{\widetilde{w}(x,0)t^{-\alpha}}{\Gamma(1-\alpha)},
\end{eqnarray*}
the first equation of \eqref{adj-prob-2} can be written as
\begin{eqnarray*}
   && \partial_t^\alpha \widetilde{w}(x,t)-\Delta \widetilde{w}(x,t)=-\frac{\widetilde{w}(x,0)t^{-\alpha}}{\Gamma(1-\alpha)},\qquad x\in\Omega,\quad t\in(0,T].
\end{eqnarray*}
The initial condition can be evaluated in the following way:
\begin{align*}
\lim\limits_{t\to 0^+}I_{0^+}^{1-\alpha}\widetilde{w}(x,t) & \approx I_{0^+}^{1-\alpha}\widetilde{w}(x,\tau)=\frac{1}{\Gamma(1-\alpha)}  \int_0^{\tau} \frac{\widetilde{w}(x,s)}{(\tau-s)^{\alpha}}\mathrm{d}s \\
&  \approx\frac{1}{\Gamma(1-\alpha)} \frac{\widetilde{w}(x,0)+\widetilde{w}(x,\tau)}{2} \int_0^{\tau} (\tau-s)^{-\alpha} \mathrm{d}s \\
& =\frac{\tau^{1-\alpha}}{2\Gamma(2-\alpha)}\big( \widetilde{w}(x,0)+\widetilde{w}(x,\tau) \big),
\end{align*}
where $\tau$ is the time step.
The fully discrete scheme to \eqref{adj-prob-2} reads:
\begin{eqnarray*}
  && (\widetilde{w}_h^0,v_h ) + (\widetilde{w}_h^1,v_h ) =\frac{2\Gamma(2-\alpha)}{\tau^{1-\alpha}}(u(f)(x)-g^\delta(x),v_h), \quad \forall v_h\in V_h,\\
  &&  (\widetilde{w}_h^{k+1},v_h) + \eta(\nabla \widetilde{w}_h^{k+1},\nabla v_h) = \sum\limits_{j=0}^{k-1}(b_j-b_{j+1}) (\widetilde{w}_h^{k-j},v_h) + b_k(\widetilde{w}_h^0,v_h) - \eta\frac{t_{k+1}^{-1}}{\Gamma(1-\alpha)}(\widetilde{w}_h^0,v_h),\quad \forall v_h\in V_h.
\end{eqnarray*}
Although an error analysis for the above method is missing,
our numerical results to be presented below show that such approximation is effective.

\subsection{Numerical results}

We first consider a one-dimensional problem with $\Omega=(0,1)$ and $T=1$. We fix $\mu(t)=\sin(2\pi t)$ and consider the following three source functions:

\noindent\textbf{Example 1.}~$f_1^\dag(x)=\exp(-x)\sin(2\pi x)$.

\noindent\textbf{Example 2.}~$f_2^\dag(x)=\left\{\begin{array}{ll}
                     2x, & 0\leq x<1/2, \\
                     -2x+2, & 1/2\leq x\leq 1.
                   \end{array}\right.$

\noindent\textbf{Example 3.}~$f_3^\dag(x)=\left\{\begin{array}{ll}
                     0.25, & 1/4\leq x\leq3/4, \\
                     0, & \text{otherwise}.
                   \end{array}\right.$

The spatial mesh consists of 40 equally spaced subintervals, i.e., $N=40, h=1/N$. The time discretization makes use of
the time step size $\tau=1/M$ with $M=50$.
The exact date $g(x)$ is obtained through solving the direct problem \eqref{dir-prob} with source function $f(x)$.
The noisy data $g^\delta$ is generated by adding a random perturbation, i.e.,
\begin{eqnarray*}
   && g^\delta=g+\delta_{\mathrm{rel}} {r\over\|r\|_{L^2(\Omega)}}\|g\|_{L^2(\Omega)} \quad\text{with}\quad r=2\text{rand}(\text{size(g)})-1,
\end{eqnarray*}
where $\delta_{\mathrm{rel}}>0$ is the relative noise level, $\delta:=\delta_{\mathrm{rel}}\|g\|_{L^2(\Omega)}$.
For the regularization parameter $\beta$, we choose $\beta\approx\delta^2$ although other choices
are possible.
The parameters $\varsigma_0$ and $\upsilon_0$ must be taken so that condition \eqref{cond1-AG-conv} is satisfied.
Based on the preliminary calculations, we have $\|u_h^{K_\tau}(f_1^\dag)\|_{L^2(\Omega)}/\|f_1^\dag\|_{L^2(\Omega)}\approx0.0066$, $\|u_h^{K_\tau}(f_2^\dag)\|_{L^2(\Omega)}/\|f_2^\dag\|_{L^2(\Omega)}\approx0.0748$ and $\|u_h^{K_\tau}(f_3^\dag)\|_{L^2(\Omega)}/\|f_3^\dag\|_{L^2(\Omega)}\approx0.0683$.
Thus we set accordingly $\varsigma_0=300$ and $\upsilon_0=10^{-4}$ in all the three tests.
The initial guess is: $f^0(x)\equiv 0, \rho^0(x)\equiv 0.5$.
As in \cite{Tian-Yuan-2019,Hank-Grpetsch-1998,Buccini-Donatelli-2017},
the algorithm will be stopped if the iteration $n$ reaches the prefixed maximum $N_{\max}=5000$
or the criterion:
\begin{eqnarray*}
   && \frac{\|f^{n}-f^{n-1}\|_{L^2(\Omega)}}{\|f^{n}\|_{L^2(\Omega)} }\leq 10^{-4} \quad\text{or}\quad    \|u_h^{K_\tau}(f^n)-g^\delta\|_{L^2(\Omega)}\leq 1.2\delta.
\end{eqnarray*}
It is notable that the second stopping criterion is the
classical Morozov's discrepancy principle \cite[Chapter 4]{Engl-Hanke-Beyvayer-1996}.

The proposed algorithm is validated through investigating the behavior of the reconstructed source error and
corresponding residue, defined by
\begin{eqnarray*}
   &&  e_r(f^n,f^\dag):=\frac{\|f^n-f^\dag\|_{L^2(\Omega)} }{\|f^\dag\|_{L^2(\Omega)}} \quad\text{and}\quad \mathrm{res}(f^n,g^\delta):=\|u_h^{K_\tau}(f^n)-g^\delta \|_{L^2(\Omega)}.
\end{eqnarray*}
The obtained numerical results are shown in Tables \ref{Tab1}$-$\ref{Tab3}, in which the impact on the reconstruction
quality of the noise level,  regularization parameters, and time fractional order is presented.
Theoretically, the choice of the regularization parameters should depend on the regularity of the
target solution: smoother is the solution, larger is $\beta$ and smaller is $\gamma$.
This expectation is consistent with the results in Table \ref{Tab1}$-$\ref{Tab3}.
The reconstructed sources listed in these tables are very satisfactory with expected accuracy.
Precisely, the smaller the noise level of data $g^\delta$, the more accurate the regularized solution
and the larger the number of iterations.
In addition, we observe that the recovery quality
gets worse as the order of the derivative order $\alpha$ decreases, probably due to lower regularity of
the solution for smaller $\alpha$.

\begin{table}[h!]\small
\setlength{\belowcaptionskip}{0.2cm}
\centering
\caption{Numerical results of the reconstruction for $f^\dag=f_1^\dag$}\label{Tab1}\setlength{\belowcaptionskip}{-0.9cm}
\begin{tabular}{cccccccc}
\hline
   &    &  \multicolumn{3}{c}{$\alpha=0.3$ } &  \multicolumn{3}{c}{$\alpha=0.9$ }  \\
   \cmidrule(lr){3-5} \cmidrule(lr){6-8}
$\delta_{\mathrm{rel}}$  & $(\beta,\gamma)$  & $n$  & $e_r(f^n,f^\dag)$  &  $\mathrm{res}(f^n,g^\delta)$   & $n$  & $e_r(f^n,f^\dag)$  &  $\mathrm{res}(f^n,g^\delta)$ \\\hline
$2\%$  &  $(5,5)\times 10^{-8}$   & $57$  &  $0.0448$  &  $3.0487e$--$04$  & $32$  &  $0.0413$  &  $2.9186e$--$04$ \\
   &  $(5,10)\times 10^{-8}$  & $59$  &  $0.0567$  &  $3.0490e$--$04$  & $32$  &  $0.0484$  &  $2.9336e$--$04$ \\
   &  $(10,5)\times 10^{-8}$  & $57$  &  $0.0449$  &  $3.0517e$--$04$  & $32$  &  $0.0413$  &  $2.9203e$--$04$ \\
   &  $(10,10)\times 10^{-8}$   & $59$  &  $0.0567$  &  $3.0524e$--$04$  & $32$  &  $0.0484$  &  $2.9355e$--$04$ \\\hline
$1\%$  &  $(2,2)\times 10^{-8}$  & $82$  &  $0.0246$  &  $1.5314e$--$04$  & $49$  &  $0.0203$  &  $1.4582e$--$04$ \\
 &  $(2,5)\times 10^{-8}$  & $88$  &  $0.0423$  &  $1.5256e$--$04$  & $48$  &  $0.0259$  &  $1.4730e$--$04$ \\
 &  $(5,2)\times 10^{-8}$  & $82$  &  $0.0247$  &  $1.5333e$--$04$  & $48$  &  $0.0204$  &  $1.4592e$--$04$ \\
 &  $(5,5)\times 10^{-8}$  & $88$  &  $0.0423$  &  $1.5276e$--$04$  & $48$  &  $0.0259$  &  $1.4741e$--$04$ \\\hline
$0.5\%$  &  $(1,1)\times 10^{-8}$   & $111$  &  $0.0196$  &  $7.6418e$--$05$  & $60$  &  $0.0200$  &  $7.2923e$--$05$ \\
&  $(1,5)\times 10^{-8}$  & $132$  &  $0.0505$  &  $7.6679e$--$05$  & $62$  &  $0.0276$  &  $7.2992e$--$05$ \\
&  $(5,1)\times 10^{-8}$  & $111$  &  $0.0196$  &  $7.6659e$--$05$  & $60$  &  $0.0200$  &  $7.3064e$--$05$ \\
&  $(5,5)\times 10^{-8}$  & $134$  &  $0.0512$  &  $7.6609e$--$05$  & $62$  &  $0.0276$  &  $7.3181e$--$05$ \\\hline
\end{tabular}
\end{table}

\begin{table}[h!]\small
\setlength{\belowcaptionskip}{0.2cm}
\centering
\caption{Numerical results of the reconstruction for $f^\dag=f_2^\dag$}\label{Tab2}\setlength{\belowcaptionskip}{-0.9cm}
\begin{tabular}{cccccccc}
\hline
   &    &  \multicolumn{3}{c}{$\alpha=0.3$ } &  \multicolumn{3}{c}{$\alpha=0.9$ }  \\
   \cmidrule(lr){3-5} \cmidrule(lr){6-8}
$\delta_{\mathrm{rel}}$  & $(\beta,\gamma)$  & $n$  & $e_r(f^n,f^\dag)$  &  $\mathrm{res}(f^n,g^\delta)$   & $n$  & $e_r(f^n,f^\dag)$  &  $\mathrm{res}(f^n,g^\delta)$ \\\hline
$2\%$  &  $(1,1)\times 10^{-7}$   & $13$  &  $0.0936$  &  $1.1177e$--$03$   & $14$  &  $0.0782$  &  $9.6530e$--$04$ \\
   &  $(1,2)\times 10^{-7}$     & $14$  &  $0.0975$  &  $1.1293e$--$03$     & $14$  &  $0.0825$  &  $9.6707e$--$04$ \\
   &  $(2,1)\times 10^{-7}$     & $13$  &  $0.0936$  &  $1.1109e$--$03$     & $14$  &  $0.0782$  &  $9.6388e$--$04$ \\
   &  $(2,2)\times 10^{-7}$     & $14$  &  $0.0975$  &  $1.1297e$--$03$     & $14$  &  $0.0825$  &  $9.6566e$--$04$ \\\hline
$1\%$  &  $(1,1)\times 10^{-7}$    & $29$  &  $0.0724$  &  $5.9592e$--$04$   & $17$  &  $0.0697$  &  $5.0506e$--$04$ \\
 &  $(1,2)\times 10^{-7}$     & $32$  &  $0.0783$  &  $5.9752e$--$04$     & $17$  &  $0.0746$  &  $5.1410e$--$04$ \\
 &  $(2,1)\times 10^{-7}$     & $29$  &  $0.0724$  &  $5.9605e$--$04$     & $17$  &  $0.0697$  &  $5.0497e$--$04$ \\
 &  $(2,2)\times 10^{-7}$     & $32$  &  $0.0783$  &  $5.9767e$--$04$     & $17$  &  $0.0746$  &  $5.1400e$--$04$ \\\hline
$0.5\%$  &  $(1,1)\times 10^{-8}$    & $52$  &  $0.0490$  &  $2.9706e$--$04$  & $30$  &  $0.0454$  &  $2.5299e$--$04$ \\
&  $(1,5)\times 10^{-8}$     & $55$  &  $0.0552$  &  $2.9805e$--$04$     & $31$  &  $0.0485$  &  $2.5322e$--$04$ \\
&  $(5,1)\times 10^{-8}$     & $52$  &  $0.0491$  &  $2.9717e$--$04$     & $30$  &  $0.0454$  &  $2.5307e$--$04$ \\
&  $(5,5)\times 10^{-8}$     & $55$  &  $0.0552$  &  $2.9816e$--$04$     & $31$  &  $0.0485$  &  $2.5330e$--$04$ \\\hline
\end{tabular}
\end{table}

\begin{table}[h!]\small
\setlength{\belowcaptionskip}{0.2cm}
\centering
\caption{Numerical results of the reconstruction for $f^\dag=f_3^\dag$}\label{Tab3}\setlength{\belowcaptionskip}{-0.9cm}
\begin{tabular}{cccccccc}
\hline
   &    &  \multicolumn{3}{c}{$\alpha=0.3$ } &  \multicolumn{3}{c}{$\alpha=0.9$ }  \\
   \cmidrule(lr){3-5} \cmidrule(lr){6-8}
$\delta_{\mathrm{rel}}$  & $(\beta,\gamma)$  & $n$  & $e_r(f^n,f^\dag)$  &  $\mathrm{res}(f^n,g^\delta)$   & $n$  & $e_r(f^n,f^\dag)$  &  $\mathrm{res}(f^n,g^\delta)$ \\\hline
$1\%$  &  $(5,5)\times 10^{-8}$    & $116$  &  $0.2619$  &  $1.7268e$--$04$   & $78$  &  $0.2446$  &  $1.4823e$--$04$ \\
   &  $(5,10)\times 10^{-8}$     & $147$  &  $0.2512$  &  $1.7267e$--$04$     & $86$  &  $0.2372$  &  $1.4842e$--$04$ \\
   &  $(10,5)\times 10^{-8}$     & $117$  &  $0.2616$  &  $1.7250e$--$04$     & $78$  &  $0.2446$  &  $1.4830e$--$04$ \\
   &  $(10,10)\times 10^{-8}$     & $148$  &  $0.2510$  &  $1.7261e$--$04$     & $86$  &  $0.2373$  &  $1.4849e$--$04$ \\\hline
$0.5\%$  &  $(5,5)\times 10^{-9}$   & $266$  &  $0.2429$  &  $8.6367e$--$05$   & $166$  &  $0.2325$  &  $7.4150e$--$05$ \\
 &  $(5,50)\times 10^{-9}$     & $367$  &  $0.2180$  &  $8.6328e$--$05$     & $191$  &  $0.2178$  &  $7.4196e$--$05$ \\
 &  $(50,5)\times 10^{-9}$     & $268$  &  $0.2428$  &  $8.6278e$--$05$     & $166$  &  $0.2325$  &  $7.4240e$--$05$ \\
 &  $(50,50)\times 10^{-9}$     & $370$  &  $0.2179$  &  $8.6333e$--$05$     & $192$  &  $0.2177$  &  $7.4173e$--$05$ \\\hline
$0.1\%$  &  $(1,1)\times 10^{-9}$    & $2278$  &  $0.1781$  &  $1.7272e$--$05$   & $1468$  &  $0.1678$  &  $1.4851e$--$05$ \\
&  $(1,10)\times 10^{-9}$    & $3378$  &  $0.1453$  &  $1.7274e$--$05$     & $1924$  &  $0.1441$  &  $1.4846e$--$05$ \\
&  $(10,1)\times 10^{-9}$     & $2301$  &  $0.1780$  &  $1.7274e$--$05$     & $1477$  &  $0.1677$  &  $1.4849e$--$05$ \\
&  $(10,10)\times 10^{-9}$     & $3398$  &  $0.1455$  &  $1.7275e$--$05$     & $1934$  &  $0.1442$  &  $1.4848e$--$05$ \\\hline
\end{tabular}
\end{table}

\begin{figure}[h!]
 \centerline{\scalebox{0.38}{\includegraphics{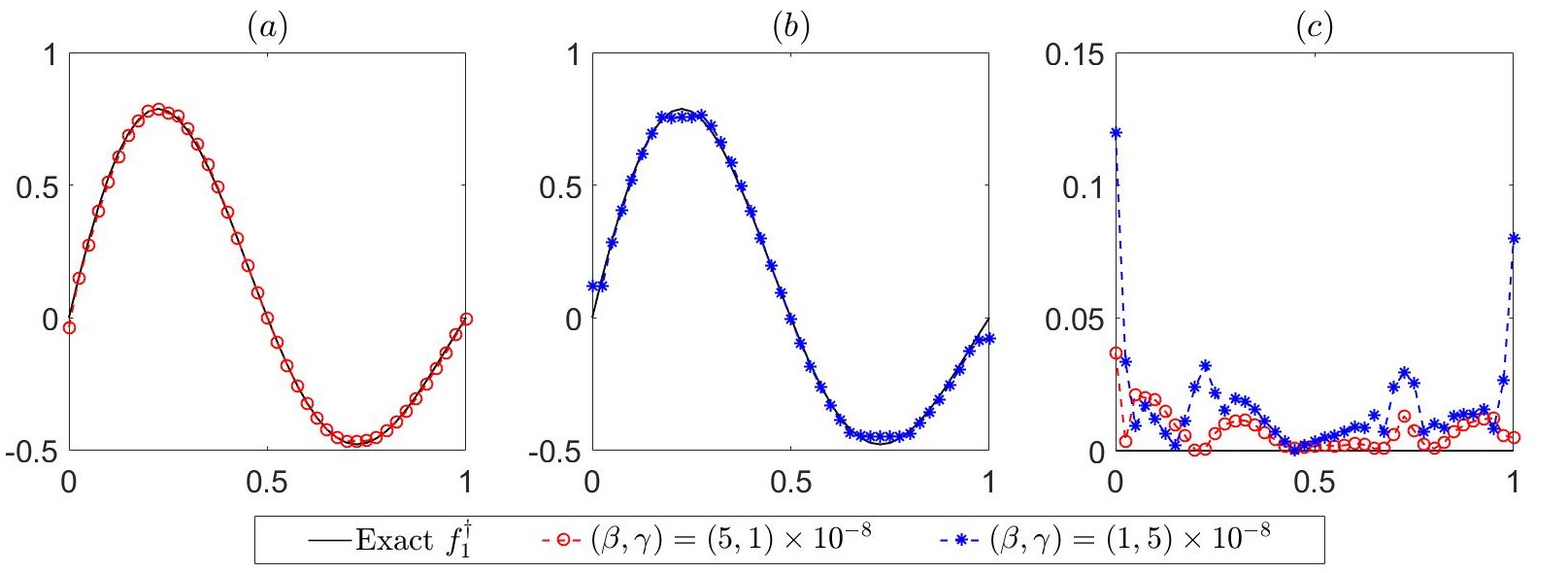}}}
  \centering\caption{\small The computed source functions for $f_1^\dag$ with $\alpha=0.9$ and $\delta_{\mathrm{rel}}=0.5\%$. $(a)$ and $(b)$: the regularized solutions;
  $(c)$: the absolute error $|f^n-f_1^\dag|$.} \label{Fig_1D_f_1}
\end{figure}

\begin{figure}[h!]
 \centerline{\scalebox{0.38}{\includegraphics{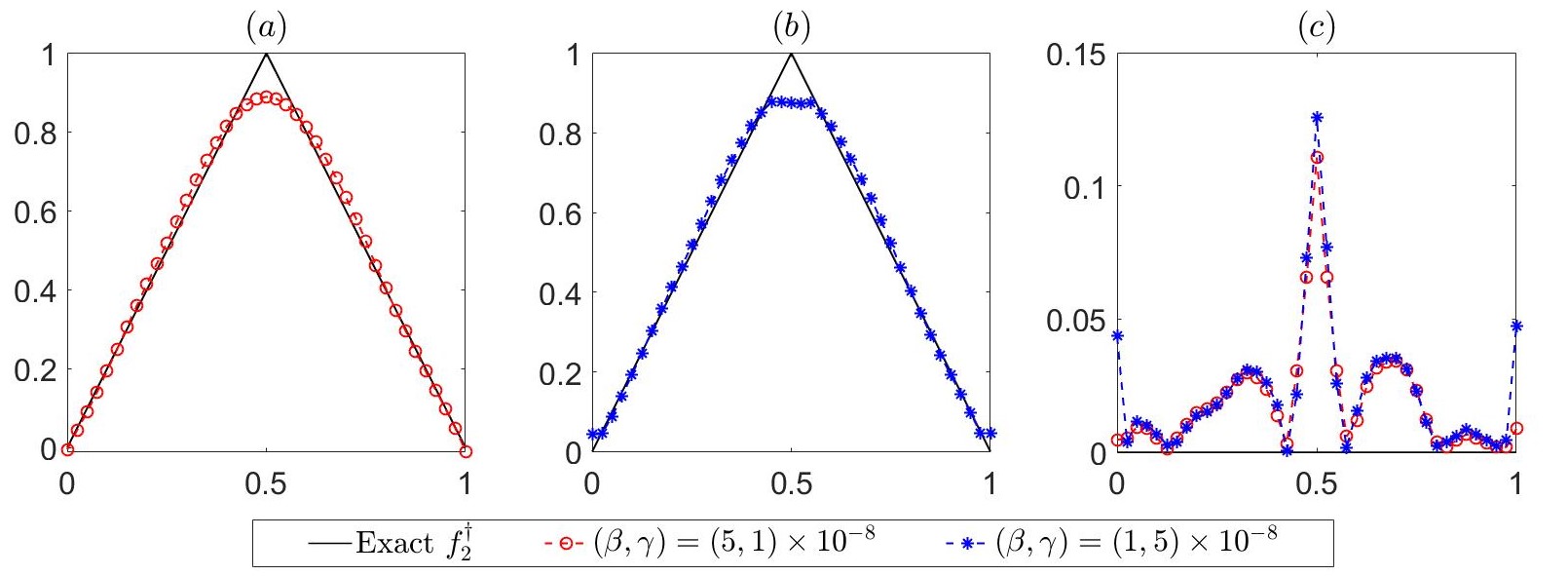}}}
  \centering\caption{\small The computed source functions for $f_2^\dag$ with $\alpha=0.9$ and $\delta_{\mathrm{rel}}=0.5\%$. $(a)$ and $(b)$: the regularized solutions; $(c)$: the absolute error $|f^n-f_2^\dag|$.} \label{Fig_1D_f_2}
\end{figure}

\begin{figure}[h!]
 \centerline{\scalebox{0.38}{\includegraphics{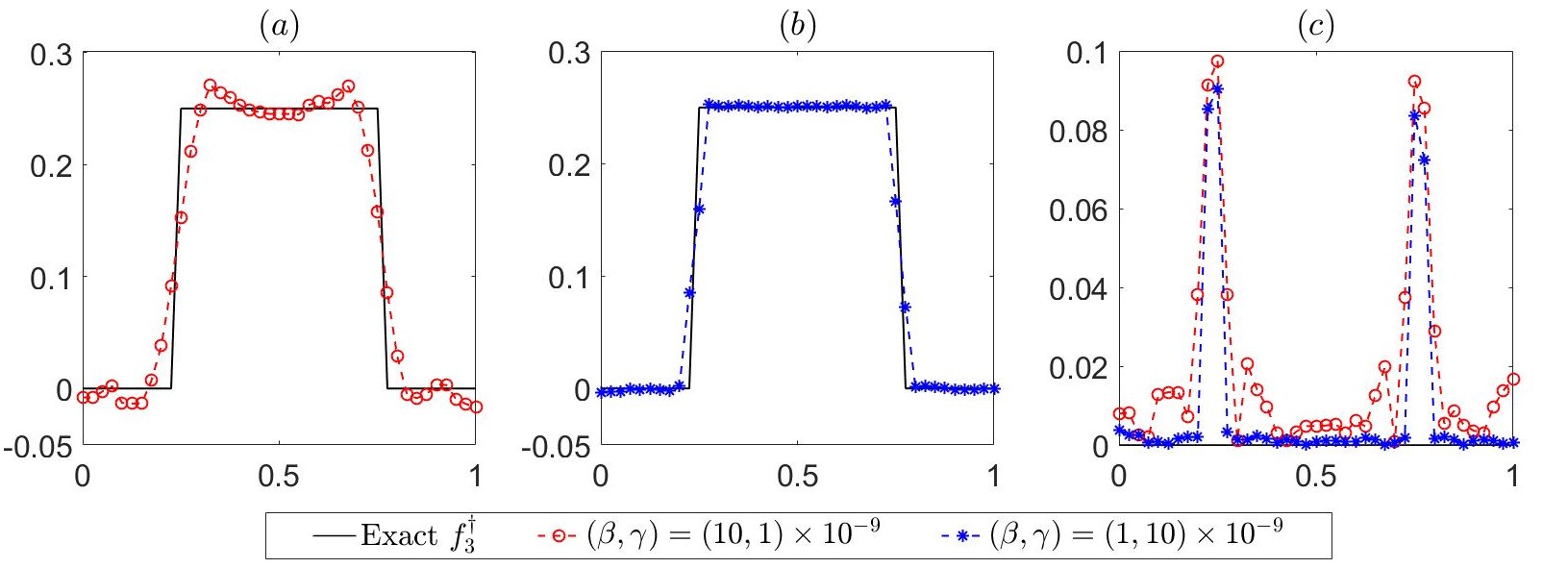}}}
  \centering\caption{\small The computed source functions for $f_3^\dag$ with $\alpha=0.9$ and $\delta_{\mathrm{rel}}=0.1\%$. $(a)$ and $(b)$: the regularized solutions; $(c)$: the absolute error $|f^n-f_3^\dag|$.} \label{Fig_1D_f_3}
\end{figure}

The computed sources and corresponding absolute errors are plotted in
Figure \ref{Fig_1D_f_1}$-$\ref{Fig_1D_f_3}. More precisely, Figure \ref{Fig_1D_f_1}(a)(b)
compares the computed results between two choices of the regularization parameters for the smooth solution
$f_1^\dag$. Much more accurate solution is obtained by using the parameter pair
$(\beta,\gamma)=(5,1)\times 10^{-8}$ than $(\beta,\gamma)=(1,5)\times 10^{-8}$.

The computed result for the piecewisely smooth solution $f_2^\dag$ is shown in
Figure \ref{Fig_1D_f_2}. For this $H^1$-regularity only solution,
the choice of the regularization parameters makes no significant difference on the accuracy, but
the parameter pair $(\beta,\gamma)=(5,1)\times 10^{-8}$ gives still slightly better result than
$(\beta,\gamma)=(1,5)\times 10^{-8}$.

Finally for the discontinuous solution $f_3^\dag$, the result presented in Figure \ref{Fig_1D_f_3}
clearly demonstrates
that the TV regularization term has the effect of stabilizing the discontinuous solution.
In particular, it eliminates high frequency oscillations and allows better recovering of the solution near
the discontinuous points.

Now we consider a two-dimensional problem.
The space domain $\Omega=(0,1)^2$ is partitioned into $2\times 40^2$ equal triangles
and the time finite difference discretization in the domain $[0,1]$ uses $50$ equidistant grid points.
We fix $\alpha=0.9$ and $\mu(t)=1$. The exact source to recovery is the discontinuous piecewise constant function as follows:

\noindent\textbf{Example 4.}~$f_4^\dag(x_1,x_2)=\left\{\begin{array}{ll}
                     0.25, & (x_1-0.5)^2+(x_2-0.5)^2\leq (0.25)^2, \\
                     0, & \text{otherwise}.
                   \end{array}\right.$

\begin{figure}[h!]
 \centerline{\scalebox{0.38}{\includegraphics{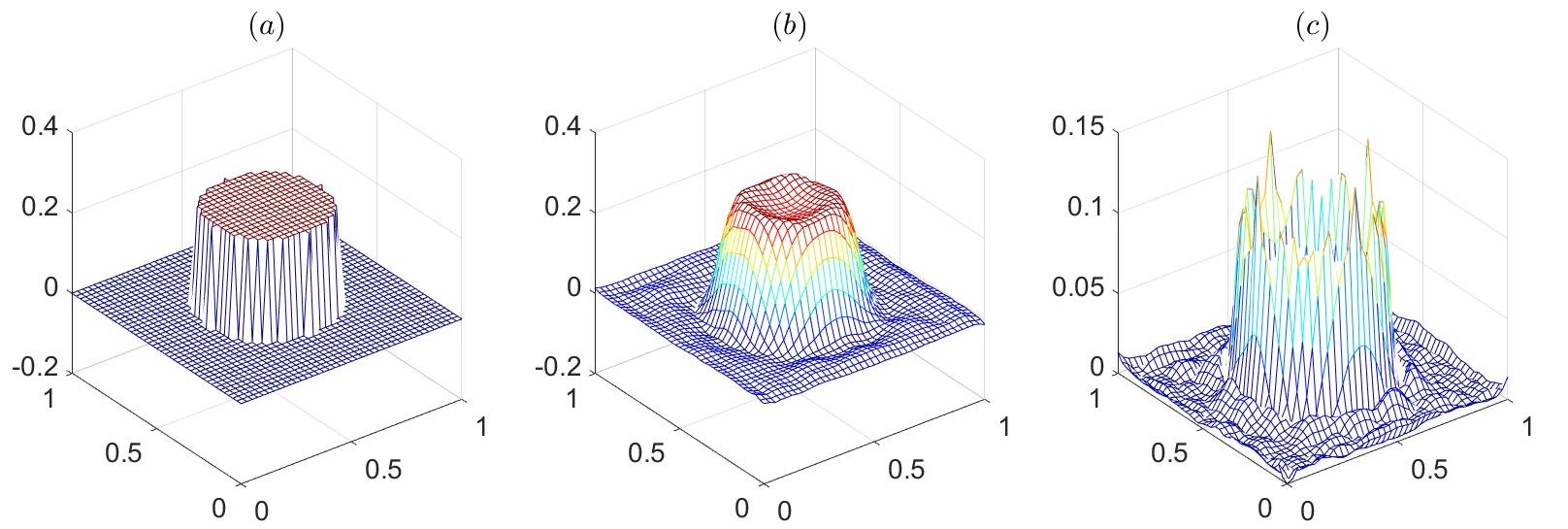}}}
 \centerline{\scalebox{0.38}{\includegraphics{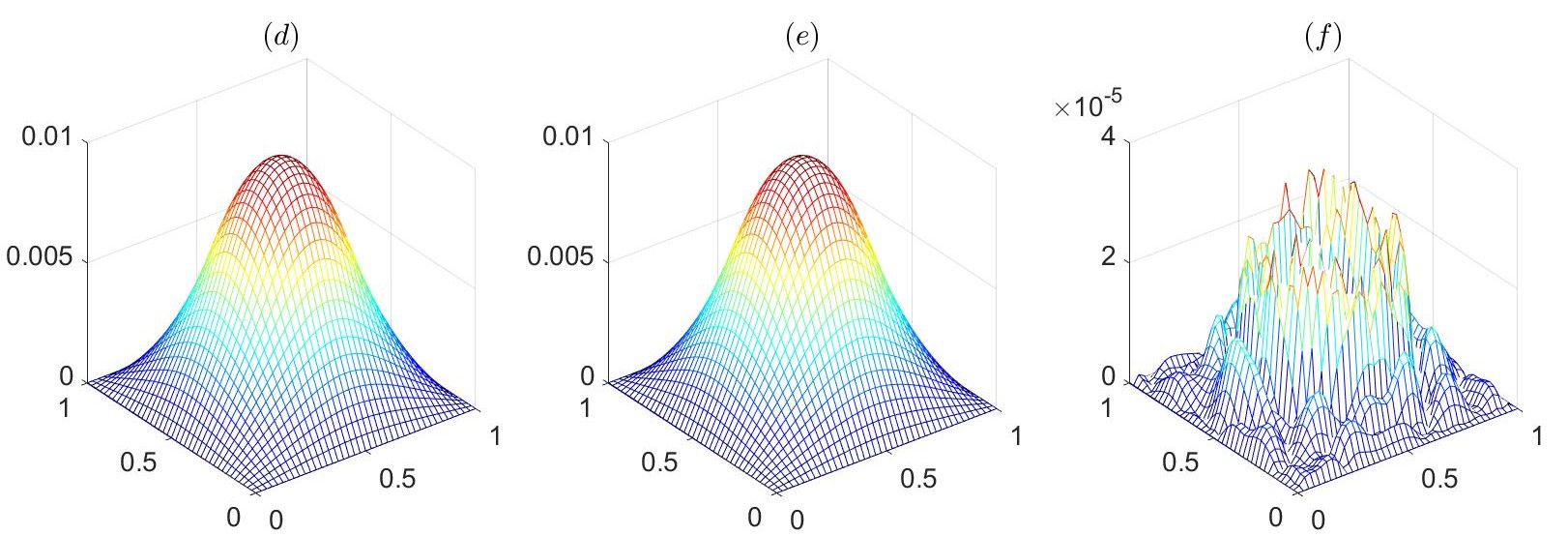}}}
  \centering\caption{\small The computed source functions for $f_4^\dag$ with $\alpha=0.9$ and $\delta_{\mathrm{rel}}=0.1\%$. $(a)$: exact solution $f_4^\dag$; $(b)$: regularized solution $f^n$;
  $(c)$: $|f^n-f_4^\dag|$;
  $(d)$: reference data $u_h^{K_\tau}(f_4^\dag)$; $(e)$: computed data $u_h^{K_\tau}(f^n)$; $(f)$: $|u_h^{K_\tau}(f^n)-u_h^{K_\tau}(f_4^\dag)|$.} \label{Fig_2D}
\end{figure}

\noindent For the parameters in Algorithm \ref{alg-acc-pd}, we set $\varsigma_0=500$, $\upsilon_0=10^{-4}$,
$(\beta,\gamma)=(1,10)\times 10^{-10}$.
The initial guess is set to $f^0(x_1,x_2)\equiv 0, \rho^0(x_1,x_2)\equiv 0.5$.
For ill-posed problems, if the solution to restore has a low regularity,
the convergence can be arbitrarily slow. This can make the algorithm very expensive,
particularly for high-dimensional problems.
In order to save computation, we set the maximum iteration number $N_{\max}=1000$.
Figure \ref{Fig_2D} presents the computed result at $n=1000$.
The set of figures in the first row represents respectively the exact source function, the
reconstructed source, and their discrepancy. The figures in the second row represent
the data, i.e., the final state solutions $u_h^{K_\tau}$, associated respectively to the exact source
$f_4^\dag)$, computed source $f^n$, and their absolute errors.
We observe satisfactory reconstruction with desirable accuracy.
The error peak appears in the neighborhood of discontinuities with about 10\% error, while the error
in the smooth region is below 2\%.
Precise values hidden behind these figures are:
$\delta=4.2832e$--$06$ for the noise level, $e_r(f^n,f_4^\dag)=0.1484$ for the source error,
and $\mathrm{res}(f^n,g^\delta)=1.6091e$--$05$ for the data error.

\section{Concluding remarks}\label{sec:6}
In this paper, we have studied
the inverse problem of recovering a source term in the time-fractional diffusion equation
of $\alpha$ order.
Unlike the existing work for this type of problems,
we proposed a regularized model with $L^2$-TV regularization, which is beneficial for reconstructing discontinuous or piecewise constant solutions.
By applying the standard Galerkin based piecewise linear finite element method in space
and the popular $2-\alpha$ order finite difference scheme in time, a fully discrete problem
for the regularized model was derived.
Regarding the theoretical aspect, we first established the convergence order of the discrete
problem to the continuous direct problem. Then
the convergence rate of the discrete regularized solution to the target solution was derived.
The convergence of the regularized solution with respect to the noise level was also provided.
Finally, in order to efficiently implement the discrete regularized model,
we proposed a primal-dual iterative algorithm based on an equivalent saddle-point reformulation
of the regularized model.
Several numerical examples are given to support the theoretical results and
verify the efficiency of the proposed method.
There remains some interesting questions for the future work.
For example, how to choose the optimal regularization parameters $\beta$ and $\gamma$ in practical calculation,
and derive the convergence rate of the iterative solutions to the exact solution, etc.

\end{document}